\documentclass[journal]{IEEEtran} 
\IEEEoverridecommandlockouts    

\usepackage[dvipsnames]{xcolor}

\usepackage{subfigure} \usepackage{mdframed} 

\usepackage{amsmath,amssymb} \usepackage{graphicx,epsfig,framed}
\usepackage{mathptmx,times} 
\usepackage{epstopdf,textcase} 
\usepackage{cite,url} \usepackage[latin1]{inputenc}
\usepackage{alltt} 
\usepackage{enumerate,siunitx} \usepackage{multirow} %
\usepackage{scalerel,stackengine} \stackMath
\newtheorem{theorem}{Theorem}[section]
 
\newtheorem{proposition}[theorem]{Proposition}

\newtheorem{remark}[theorem]{Remark} 
\newtheorem{assumption}{Assumption}

\newcommand{\BE}{\begin{equation}} \newcommand{\BEQ}[1]{\BE\mathlabel{#1}} %
\newcommand{\EEQ}{\end{equation}} \newcommand{\rfb}[1]{\mbox{\rm
(\ref{#1})}\ifx\undefined\stillediting\else:\fbox{$#1$}\fi}

\newcommand{\rvline}{\hspace*{-\arraycolsep}\vline \hspace*{-\arraycolsep}} %
\newcommand{\rline}  {{\mathbb R}} 
\newcommand{\AAA}  {{\mathbf A}}

 \newcommand{\HHH}  {{\mathbf H}}

\newcommand{\PPP}  {{\mathbf P}} 
 
\newcommand{\xxx}  {{\mathbf x}} \newcommand{\zzz}  {{\mathbf z}} %
    
  \newcommand{\Lscr}
{{\cal L}}  
\newcommand{\Sscr} {{\cal S}}    
\newcommand{\mm}     {{\hbox{\hskip 0.5pt}}} \newcommand{\m}      {{\hbox{\hskip
1pt}}}  \newcommand{\nm}    
{{\hbox{\hskip -3pt}}} \newcommand{\bluff}  {{\hbox{\raise 15pt \hbox{\mm}}}}
\newcommand{\sbluff} {{\hbox{\raise  9pt \hbox{\mm}}}} %
\renewcommand{\L}    {{\Lambda}} 
\renewcommand{\o}    {{\omega}} \newcommand{\e}      {{\varepsilon}}
 
\newcommand{\dd}     {{\rm d\hbox{\hskip 0.5pt}}} \newcommand{\FORALL}
{{\hbox{$\hskip 11mm \forall \;$}}} \newcommand{\rarrow} {\mathop{\rightarrow}}
 \newcommand{\Om}     {\Omega} %
  
\newcommand\reallywidehat[1]{%
\savestack{\tmpbox}{\stretchto{%
\scaleto{%
	\scalerel*[\widthof{\ensuremath{#1}}]{\kern-.6pt\bigwedge\kern-.6pt}%
	{\rule[-\textheight/2]{1ex}{\textheight}}
}{\textheight}%
}{0.5ex}}%
\stackon[1pt]{#1}{\tmpbox}%
} 

\let\oldlabel=\label

\renewcommand{\label}[1]{\leavevmode\smash{\raise 10pt\llap
{\fbox{\scriptsize#1}}}\oldlabel{#1}} \newcommand{\mathlabel}[1]{\smash{\raise
9pt\llap {\scriptsize(#1)}}\label{#1}}

\renewcommand{\label}[1]{\oldlabel{#1}} \renewcommand{\mathlabel}[1]{\label{#1}}

\newcommand{\bbm}[1]{\left[\begin{matrix} #1 \end{matrix}\right]}
\newcommand{\sbm}[1]{\left[\begin{smallmatrix} #1 \end{smallmatrix}\right]}

\begin{document}

\title{The equilibrium points and stability of grid-connected synchronverters}

\author{ \vskip 1em Pietro Lorenzetti, Zeev Kustanovich, Shivprasad Shivratri
and George Weiss \thanks{This research was supported by the Israeli Ministry of
Infrastructure, Energy and Water, grant numbers 217-11-037 and 219-11-128. P.~Lorenzetti is a
team member in the ITN network ConFlex, funded by the European Union's Horizon
2020 research and innovation program under the Marie Sklodowska-Curie grant
agreement no. 765579.

P. Lorenzetti, S. Shivratri and G. Weiss are with the School of Electrical
Engineering, Tel Aviv University, Ramat Aviv, Israel (e-mail:
shivprasadshivratri8793@gmail.com, plorenzetti@tauex.tau.ac.il,
gweiss@tauex.tau.ac.il)

Z. Kustanovich is with the Israel Electricity Company, North District (e-mail:
kustanz875@gmail.com, phone: +972-523995779)}} \maketitle


\begin{abstract} Virtual synchronous machines are
inverters with a control algorithm that causes them to behave towards the power
grid like synchronous generators. A popular way to realize such inverters are
synchronverters. Their control algorithm has evolved over time, but all the
different formulations in the literature share the same ``basic control
algorithm''. We investigate the equilibrium points and the stability of a
synchronverter described by this basic algorithm, when connected to an infinite
bus. We formulate a fifth order model for a grid-connected synchronverter and
derive a necessary and sufficient condition for the existence of equilibrium
points. We show that the set of equilibrium points with positive field current
is a two-dimensional manifold that can be parametrized by the corresponding pair
$(P,Q)$, where $P$ is the active power and $Q$ is the reactive power. This
parametrization has several surprizing geometric properties, for instance, the
prime mover torque, the power angle and the field current can be seen directly
as distances or angles in the $(P,Q)$ plane. In addition, the stable equilibrium
points correspond to a subset of a certain angular sector in the $(P,Q)$ plane.
Thus, we can predict the stable operating range of a synchronverter from its
parameters and from the grid voltage and frequency. Our stability result is
based on the intrinsic two time scales property of the system, using tools from
singular perturbation theory. We illustrate our theoretical results with two
numerical examples. \end{abstract}

\begin{IEEEkeywords} Virtual synchronous machine, frequency droop, voltage
droop, inverter, synchronverter, Park transformation, saturating integrator,
singular perturbation method. \end{IEEEkeywords}



\section{Introduction} \label{sec1}  

\IEEEPARstart{M}{o}st distributed generators are connected to the utility grid
via inverters that rely on various control algorithms to maintain synchronism. 
They usually offer no inertia, and behave as controlled current sources that
produce fluctuating power. Numerous researchers are investigating how the future
power grids should be controlled when inverters become dominant, offering
competing control algorithms for grid-forming converters, see for instance the
recent study \cite{Tayyebi2020a}. One of the proposed approaches is to emulate
the behavior of synchronous generators (SG), so that an inverter-based grid
behaves like one based on SG, see for instance
\cite{Beck2007,DongChen:17,Driesen2008,
Mandrile2021,Mo2017, Roldan-Perez2019,Tayyebi2020a,Wu2016,Zhong2009}. This has many advantages, such as backward compatibility with
the current grid, well known black start and fault ride-through procedures, and
well tested primary and secondary frequency and voltage support algorithms.
Following \cite{Beck2007}, inverters that behave towards the utility grid like
synchronous machines are called {\em virtual synchronous machines} (VSM).

\begin{figure}[h] \centering 
\includegraphics[width=1\linewidth]{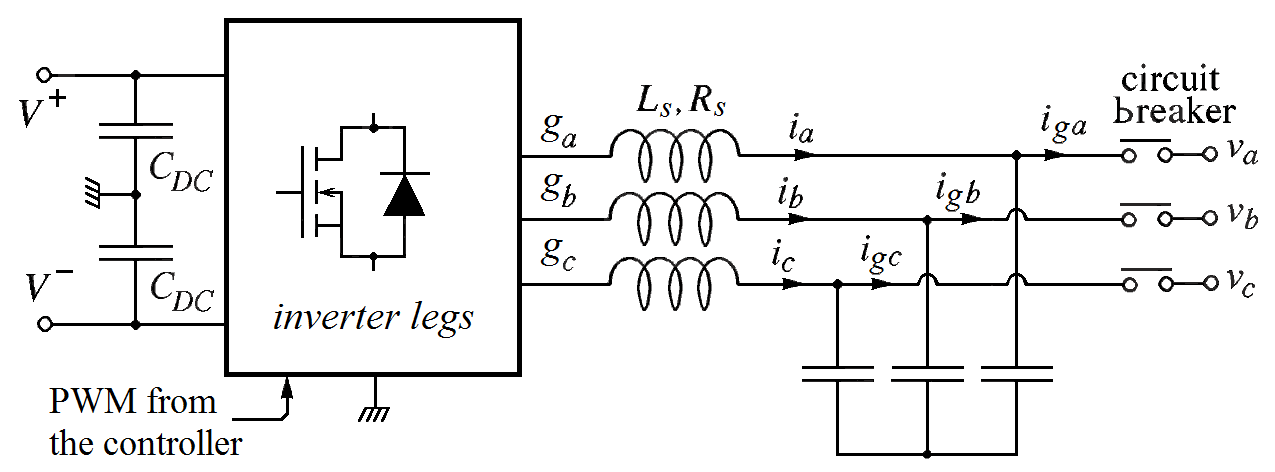}
\caption{An inverter with an LC filter receiving DC
	voltages $V^+$, $V^-$ and connected to the grid voltages $v_a,\m v_b,\m v_c$.
	The control algorithm receives measurements of $v_a,\m v_b,\m v_c$ and of the
	grid currents $i_{ga},\m i_{gb},\m i_{gc}$.} \label{fig:InverterLC}
\end{figure}

One particular type of VSM are the {\em synchronverters}, introduced
in \cite{Zhong2009,Zhong2011}. This type of inverter
has attracted considerable attention, see for instance
 \cite{Alipoor2013,Aouini2015,Blau2018,Dong2016,Shuai2020,
Roldan-Perez2019,Zhong2013,Zhong2016,Zhong2014}, and the recent survey
\cite{Vasudevan2020}. The hardware of a synchronverter is similar to that of a
conventional three phase inverter (with any number of DC
levels, most commonly 3), the novelty is in the control algorithm. The only
hardware difference is that some fast acting energy storage (typically,
capacitors) is required on the DC bus, to provide the energy pulses (both
positive and negative) needed for the emulation of rotor inertia. We base our
modelling on the simplified circuit diagram of a grid-connected inverter
in Fig.~\ref{fig:InverterLC}. Even though the synchronverter
control scheme has evolved over time, all the formulations present in
the literature share the same ``basic control algorithm''. We base our modelling
on this basic algorithm (see Fig.~\ref{fig:block} in Sect. \ref{sec2} for more
details).

The stability of a SG (or VSM) connected to an infinite bus is a well motivated
classical problem in the study of power systems. For instance, in
\cite[Sect.~12.3]{Kundur1994} and \cite{Liu2019} we can find
the stability analysis of a linearized second order approximation of this
system. The same problem, with a more complex SG model, has been addressed in
\cite{Sauer1997}. In the last decade, motivated by the
growing interest in VSM-based grids, similar stability studies have been
performed for more complex models of grid-connected VSM. In
\cite{Barabanov2017}, \cite{Natarajan2018} a fourth order model is formulated
for a grid-connected VSM, and conditions on the parameters ensuring almost
global asymptotic stability (aGAS) are derived, and in
\cite{Mandrile2020} a novel technique for state-space modeling of grid-connected
converters is presented, with local stability evaluation via eigenvalues. In
\cite{Tayyebi2020b} a VSM model is developed, which contains a DC side that
interacts with the AC side in an ingenious way, leading to aGAS of the VSM
connected to an infinite bus. In the context of microgrids, stability results
are derived in \cite{Bretas2003,Schiffer2014,Zhang2016}, and the importance of
accurate modeling has been discussed, among others, in \cite{Vorobev2018}
and in the recent review \cite{Ojo2020}. The paper 
\cite{Milani2018} presents an interesting equilibrium point analysis for 
microgrids interfaced via solid-state transformers. This analysis is then employed
to develop a power sharing algorithm between the inverters.

This paper investigates the local asymptotic stability of a VSM functioning
according to the basic synchronverter algorithm, when connected to a powerful
grid modelled as an infinite bus. For this purpose, we formulate a fifth order
grid-connected synchronverter model. This model is an
extension of the fourth order model developed and analyzed in
\cite{Barabanov2017, Natarajan2018}, where the rotor (or field) current was
assumed to be constant (thus ignoring the reactive power control loop). Using
advanced mathematical methods, different sufficient conditions for almost global
asymptotic stability of the fourth order model were derived in
\cite{Barabanov2017, Natarajan2018}. Here we include the field current as the
fifth state variable and we investigate the stability of the equilibrium points
of the resulting fifth order system.

We derive a novel geometric representation of the fourth and fifth order models'
equilibrium points. We use extensively the mapping of equilibrium points into
the {\em power plane}, where the coordinates are $P$ (the active power) and $Q$
(the reactive power). (In the language of differential geometry, the manifold of
equilibrium points with positive field current is diffeomorphic to the power
plane.) We show that, for a fixed prime mover torque and for field current
values in a ``reasonable'' range, the image of the fourth order model
equilibrium points in the power plane moves on a circle. The radius of this
circle depends on the prime mover torque at the equilibrium. In the same
geometric representation, we identify a stability sector for the fifth order
model equilibrium points. This sector allows to determine a priori if certain
reference values of active and reactive power will generate stable (or unstable)
fifth order model equilibrium points.

The paper is organized as follows. In Sect. \ref{sec2} we recover the fourth
order grid-connected synchronverter model from
\cite{Natarajan2018,Natarajan2017}, and we extend it to a fifth order one,
adding the field current to the state vector. In Sect. \ref{sec3} the
equilibrium points of the fourth order model are studied and the novel geometric
representation is introduced. In Sect. \ref{sec4} we study the equilibrium
points of the fifth order model and their representation in the power plane.
In Sect. \ref{sec5} we use results from Sect. \ref{sec3} and \ref{sec4} to find
a sufficient condition ensuring the stability of the fifth order model
equilibrium points, employing singular perturbation methods developed in
\cite{Lorenzetti2020}. Based on this result, we characterize the power plane
region corresponding to stable fifth order model equilibrium points. Finally, in
Sect. \ref{sec6} we use two numerical examples to illustrate our novel geometric
representation and our theoretical derivations.

\section{Modelling the grid-connected synchronverter} \label{sec2} 

In this section we construct the basic fifth order model of the synchronverter,
following the terminology and notation of
\cite{Natarajan2018,Natarajan2017,Zhong2011}. Note that
the paper \cite{Natarajan2017} has proposed five modifications to the
synchronverter algorithm from \cite{Zhong2011}, to improve its stability and
performance. Of these, we adopt here only the two most important ones: a
substantial increase of the effective size of the filter inductors, by using
virtual inductors, and the improved anti-windup field current controller.

\begin{figure}[h] \centering 
\includegraphics[width=0.95\linewidth]{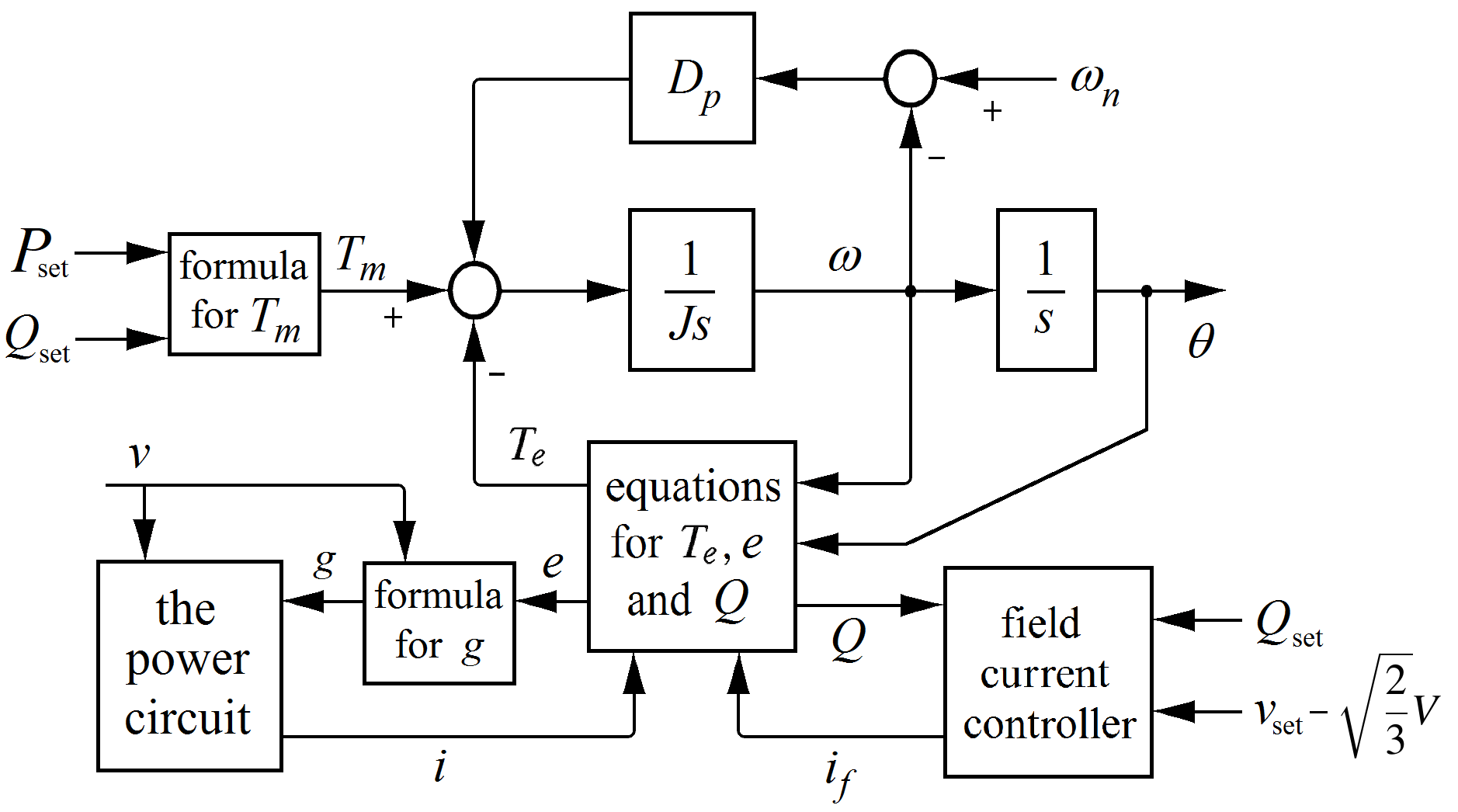} \vspace{-1mm}
\caption{The simplified block diagram of a synchronverter, adapted from
\cite[Fig.~2]{Natarajan2017}. The signals $v,g,e$ and $i$ are three
dimensional.} \label{fig:block} \vspace{-1mm} \end{figure}

Our analysis is based on a simplified model of a synchronverter, given in
Fig.~\ref{fig:block}. This model is simplified because
it does not take into account the various low-pass filters that are included to
reduce high frequency noise, and it also ignores most of the saturation blocks
included in the algorithm (see \cite{Natarajan2017,Blau2018}) (however, the
saturating integrator contained in the field current controller is considered).
We also ignore start-up procedures and various protections. Including these
features would result in a high-order model that is practically impossible to
analyze rigorously. Moreover, ignoring the aforementioned features does not
significantly alter the steady-state behaviour of the system, making our
stability analysis relevant also for higher order models.

We proceed as follows: first we recover the fourth
order model from \cite{Natarajan2017} (where the field current $i_f$ was assumed
to be a parameter). Then, we extend this model by including $i_f$ as a state
variable, obtaining a fifth order model. 

We denote by $\theta_g$ the grid angle and by $\o_g$ the grid frequency, so that
$\o_g=\dot\theta_g$. The frequency $\o_g$ is usually close to $\o_n=100\pi$\m
rad/sec. We denote by $\theta$ the synchronverter virtual rotor angle, and by
$\o$ its angular velocity, so that $\o=\dot\theta$. The difference $\delta=
\theta-\theta_g$ is called the {\em power angle}. The notation
$\widetilde\cos\m\theta$ and $\widetilde\sin\m\theta$ is defined by
\vspace{-1mm} $$ \widetilde\cos\m\theta \m=\m \left[\cos\theta\ \m
\cos\left(\theta -\frac{2\pi}{3}\right)\ \m \cos\left(\theta +
\frac{2\pi}{3}\right) \right]^\top\m,$$ \vspace{-2mm} $$ \widetilde\sin\m\theta
\m=\m \left[\sin\theta\ \m \sin\left(\theta -\frac{2\pi}{3}\right)\ \m
\sin\left(\theta + \frac{2\pi}{3} \right)\right]^\top.$$ Then the grid voltage
vector is \vspace{-1mm} \begin{equation} \label{eq:v} v \m=\m
\sqrt{\frac{2}{3}}V\widetilde\sin\m\theta_g, \end{equation} where $V$ is a
positive constant or a slowly changing signal (this is the rms value of the line
voltage).

Denote by $M_f>0$, the peak mutual inductance between the virtual rotor winding
and any one stator winding, by $i_f$ the variable {\em field current} (or rotor
current) and by $e$ the vector of electromotive forces, also called the {\em
internal synchronous voltage}. We rewrite \cite[eq.\m(4)]{Zhong2011}:
\begin{equation} \label{eq:e} e \m=\m M_f i_f\o\m\widetilde\sin\theta -
M_f\frac{\dd i_f} {\dd t} \widetilde\cos\theta \end{equation} and we note that
the variable current $i_f$ governs the amplitude of $e$. We apply the {\em Park
transformation} $$ U(\theta) = \sqrt{\frac{2}{3}} \bbm{\vspace{1mm} \cos\theta &
\cos(\theta-\frac{2\pi}{3}) & \cos(\theta+\frac{2\pi}{3}) \\ \vspace{1mm}
-\sin\theta & -\sin(\theta-\frac{2\pi}{3}) & -\sin(\theta+\frac{2\pi}{3}) \\
1/\sqrt 2 & 1/\sqrt 2 & 1/\sqrt 2} \vspace{1mm}$$ to \rfb{eq:e}. For any
$\rline^3$-valued signal $v$, the first two components of $U(\theta)v$ are
called the $dq$ coordinates of $v$, denoted by $v_d$, $v_q$. By using the
notation $m=\sqrt{3/2} \m M_f$, we represent the internal synchronous voltage
$e$ in $dq$ coordinates as: \begin{equation} \label{eq:e_qe_d} e_d \m=\m
-m\frac{\dd i_f}{\dd t}, \qquad e_q \m=\m -mi_f\o. \end{equation} The term $e_d$
can be neglected, because the rate of change of the field current is small, so
that $e_d<<e_q$. Thus, in the synchronverter algorithm from \cite{Zhong2011} the
approximation $e_d=0$ is adopted, and our analysis will follow this.
(We remark that we did simulation experiments with $e_d$
as in \rfb{eq:e_qe_d}, and the results were practically the same as for
$e_d=0$.)

Applying the Park transformation to \rfb {eq:v}, we get the $dq$ representation
of the grid voltage as \begin{equation} \label{eq:v_qv_d} v_d \m=\m
-V\sin\delta, \qquad v_q \m=\m -V\cos\delta. \end{equation}

\begin{figure}[h] \centering 
\vspace{-2mm} \includegraphics[width=0.95\linewidth]{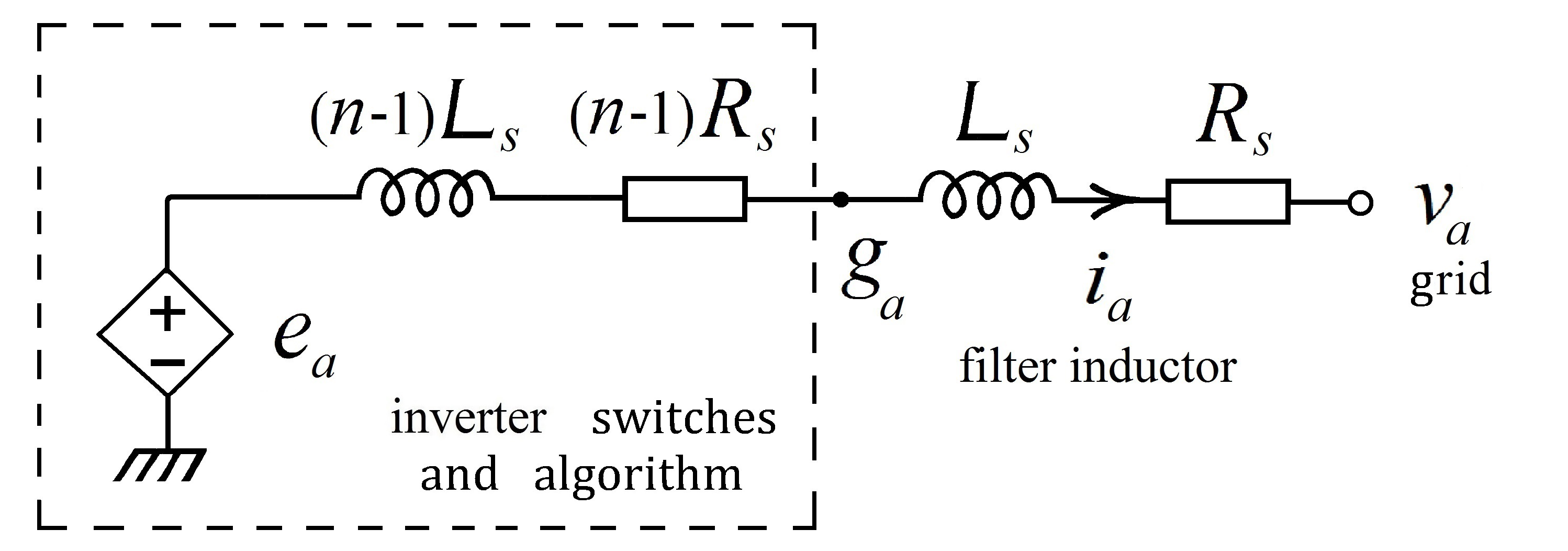}
\vspace{-2mm} \caption{A synchronverter with filter inductor $L_s$ and its
resistance $R_s$. $e_a$ is the synchronous internal voltage. The inductor and
the resistor multiplied with $(n-1)$ are virtual. Only phase $a$ is shown.
This is taken from Fig.~4 in \cite{Natarajan2017}.} \label{fig:PhaseA}
\vspace{-2mm} \end{figure}

The control algorithm computes $g=[\m g_a\ \m g_b\ \m g_c]^\top$ and sends it to
the switches in the power part (see Fig. \ref{fig:PhaseA}). In the original
synchronverter algorithm from \cite{Zhong2011}, we have $g_d=e_d$ and $g_q=e_q$,
i.e., the internal synchronous voltage $e$ from \rfb{eq:e} is sent straight to
the PWM signal generator. Here we consider the modified synchronverter equations
according to \cite{Natarajan2017}, which contain the original algorithm as a
special case, namely $n=1$. Writing \cite[eq.\m(22)]{Natarajan2017} in $dq$
coordinates, we have $$ g_d =\m \frac{(n-1)v_d + e_d}{n},\ \ \ \ \m g_q =\m
\frac{(n-1)v_q + e_q}{n} \m.$$

The current sensors are placed after the filter
capacitors, as shown in Fig.~\ref{fig:InverterLC}, to avoid the switching noise.
From these measurements, the inductor currents $i_a,\m i_b,\m i_c$ can be
estimated. (Alternatively, in some versions of synchronverters, the output
elements $L_s$ and $R_s$ are virtual, e.g., see \cite{KSW2021}, and then the currents $i_a,\m i_b,\m i_c$
are computed in the algorithm from the voltage measurements.) By applying the
Park transformation on the circuit equations corresponding to
Fig.~\ref{fig:PhaseA}, we have \begin{equation} \label{eq:id} L_s\frac{\dd
i_d}{\dd t} \m=\m -R_si_d + \o L_s i_q + g_d - v_d, \end{equation}
\begin{equation} \label{eq:iq} L_s\frac{\dd i_q}{\dd t} \m=\m -\o L_si_d - R_s
i_q + g_q - v_q.  \vspace{1mm} \end{equation} Here, $L_s$ and $R_s$ are positive
constants. Combining \rfb{eq:e_qe_d}-\rfb{eq:iq} (with $e_d=0$) and using the
notation $$R \m=\m nR_s \m,\qquad L \m=\m nL_s \m,\vspace{-3mm}$$ we get
\vspace{-1mm}\begin{equation} \label{eq:id_new} L \frac{\dd i_d}{\dd t} \m=\m
-Ri_d + \o Li_q + V\sin\delta, \vspace{-1mm} \end{equation} \begin{equation}
\label{eq:iq_new} L \frac{\dd i_q}{\dd t} \m=\m -\o Li_d - Ri_q - mi_f\o + V\cos
\delta. \end{equation}

The angular frequency evolves according to the {\em swing equation}
\vspace{-2mm} \begin{equation} \label{eq:swing} J\frac{\dd \o}{\dd t} \m=\m T_m
- T_e - D_p\o + D_p\o_n, \end{equation} where $J>0$ represents the virtual
inertia of the rotor, $T_m>0$ is the nominal active mechanical torque from the
prime mover, \BEQ{eq:Te} T_e \m=\m -mi_fi_q \end{equation} is the electric
torque computed using the measured output currents, $\o_n$ is the nominal grid
frequency and $D_p>0$ is the {\em frequency droop constant}.
We refer to \rfb{eq:T_m} for a way to choose the value of $T_m$.

We assume here that the inverter works in the linear
region of the frequency droop. The actual droop function contains dead-band and
saturation, but taking these into account would make the analysis very
complicated.

The following equation comes from the definition of $\delta$: \begin{equation}
\label{eq:delta} \frac{\dd \delta}{\dd t} \m=\m \o - \o_g \m. \end{equation}

The fourth order grid-connected synchronverter model, which considers $i_f$ as a
given parameter, can be constructed by combining the equations
\rfb{eq:id_new}-\rfb{eq:delta}. Its state vector is \begin{equation}
\label{eq:x} {\bf x} \m=\m [\m i_d \ i_q \ \o \ \delta\m]^\top \m \in \m
\mathbb{R}^4. \end{equation} We write it as a nonlinear dynamical system:
\BEQ{eq:ss4} \HHH\dot{\xxx} \m=\m \AAA(\xxx,i_f)\xxx \m + f(\xxx) \m,
\vspace{-2mm} \end{equation} where \vspace{-2mm} \begin{equation*}
\label{Lukashenko_4} \HHH \m=\m \bbm{L & 0 & 0 & 0 \\ 0 & L & 0 & 0 \\ 0 & 0 & J
& 0 \\ 0 & 0 & 0 & 1} \m,\qquad \m f(\xxx) \m=\m \bbm{V\sin \delta\\
V\cos\delta\\ T_m \m + D_p\o_n\\ \vspace{1mm} -\o_g}, \end{equation*} and $$
\AAA(\xxx,i_f) \m=\m \bbm{-R & \o L & 0 & 0 \\ -\o L & -R & -mi_f & 0 \\ 0 & m
i_f & -D_p & 0 \\ 0 & 0 & 1 & 0} \m.$$

We now derive the fifth order basic model of a synchronverter, by including
$i_f$ into the vector of the state variables. The instantaneous inverter output
reactive power is \BEQ{eq:Q} Q \m=\m v_q i_d - v_d i_q \m=\m V [i_q \sin\delta -
i_d \cos \delta] \m, \end{equation} see \cite[eq.\m(16)]{Natarajan2018}. For
convenience, we introduce \begin{equation}\label{eq:Q_tilde} \tilde Q \m=\m
Q_{\rm set} + D_q\left(\nm v_{\rm set} - \sqrt{\frac{2}{3}}V \right) \m,
\end{equation} where $v_{\rm set}$ is the desired amplitude of $v$, $D_q>0$ is
the {\em voltage droop coefficient}, $Q_{\rm set}$ is the desired reactive
power, and $V$ is as in \rfb{eq:v}. Then, the field current $i_f$ evolves
according to \begin{equation} \label{eq:Mfif} M_f \frac{\dd i_f}{\dd t} \m=\m
\frac{\tilde{Q}-Q}{K}, \end{equation} see \cite[eq.\m(15)]{Natarajan2017}, where
$K>0$ is a large constant. We want to make sure that $i_f$ stays in a reasonable
operating range $[u_{min},u_{max}]$. (We will say more about this range in Sect.
\ref{sec5}.) For this, we replace the integrator from \rfb{eq:Mfif} with a
\textit{saturating integrator} (see \cite{Lorenzetti2020}), obtaining
\cite[eq.\m(21)]{Natarajan2017}: \begin{equation} \label{eq:Mfif_sat} \frac{\dd
i_f}{\dd t} \m=\m \Sscr\left(i_f,\frac{\tilde{Q}-Q} {\tilde{K}} \right),
\end{equation} where $\tilde{K}=KM_f$. Denoting
$w=\frac{\tilde{Q}-Q}{\tilde{K}}$, the function $\Sscr$ is defined by
\vspace{-4mm} \begin{equation*} \label{eq:cal_S} \Sscr(i_f,w) \m=\m
\begin{cases} w^+ & \text{if} \quad i_f \leq u_{min}, \\ w & \text{if} \quad i_f
\in (u_{min},u_{max}), \\ w^- & \text{if} \quad i_f\geq u_{max}, \end{cases}
\end{equation*} where $w^+=\max\{w,0\}$, $w^-=\min\{w,0\}$, so that $w=w^++w^-$.
This means that as long as $i_f$ is in the range $(u_{min},u_{max})$, we have
$\frac{\dd i_f}{\dd t}=\frac{\tilde{Q}-Q}{\tilde{K}}$. However, if $i_f$ reaches
one of the end points of $[u_{min},u_{max}] $, it is not allowed to continue out
of this interval. (Note that in \rfb{eq:Mfif_sat} we use
$\tilde{K}$ in place of $K$ because, differently from \cite{Natarajan2017}, here
$i_f$ is the state, not $M_fi_f$.) Using a saturating integrator in place of a
usual one is needed in practice, and also in our stability proof in Sect.
\ref{sec5}.

The fifth order grid-connected synchronverter model can be constructed by
combining \rfb{eq:ss4}, \rfb{eq:Q}, and \rfb{eq:Mfif_sat} as: \begin{equation}
\label{eq:ss5_sat} \HHH\dot{\xxx} \m=\m \AAA(\zzz)\xxx \m + f(\xxx), \qquad
\frac{\dd i_f}{\dd t} \m=\m \Sscr\left(i_f, \frac{\tilde{Q}-Q}{\tilde{K}}
\right), \end{equation} with $\xxx$ from \rfb{eq:x}, and with the state
$\zzz=\sbm{{\bf x}\\ i_f}\in\rline^5$. Clearly, we mean that
$\AAA(\zzz)=\AAA(\xxx,i_f)$. If we ignore the saturating feature of $\Sscr$ in
\rfb{eq:Mfif_sat}, and we use \rfb{eq:Mfif} (with $Q$ from \rfb{eq:Q}) instead
of \rfb{eq:Mfif_sat} for the evolution of $i_f$ (this is true for
$i_f\in(u_{min},u_{max})$), i.e., \vspace{-1mm}
\begin{equation} \label{eq:mif} m\frac{\dd i_f}{\dd t} \m=\m k i_d \cos\delta -
k i_q \sin\delta +\frac{k}{V}\tilde{Q}, \end{equation} then we get the fifth
order non-saturated model \BEQ{eq:ss5} \tilde{\HHH}\dot{\zzz} \m=\m
\tilde{\AAA}(\zzz) \zzz \m + \tilde{f} (\zzz) \m, \end{equation} where $$
\tilde{\HHH} \m=\m \left[\begin{array}{c|c} \HHH & 0 \\ \hline 0 & m
\rule{0pt}{2ex} \end{array} \right]\m, \quad \tilde{f}(\zzz) \m=\m
\left[\begin{array}{c} f(\xxx) \vspace{1mm} \\ \frac{k}{V} \tilde{Q} \end{array}
\right], \quad k \m=\m \sqrt{\frac{3}{2}} \frac{V}{K} \m, \vspace{2mm}$$ $$
\tilde{\AAA}(\zzz)\m=\m \left[ \begin{array}{cccc|c}
\multicolumn{4}{c}{\scalebox{1}{$\AAA(\zzz)$}} \rvline \hspace{.525mm} & 0
\rule[-0.9ex]{0pt}{0pt} \\ \hline k\cos \delta & -k\sin\delta & 0 & 0 & 0
\rule{0pt}{2.6ex} \end{array} \right] \m,$$ with $\HHH$, $\AAA$, and $f$ as
defined after \rfb{eq:ss4}.

An extension of the model \rfb{eq:ss5}, to include the
effect of measurement errors, has been derived in \cite{KSW2021}. This is
needed to analyze the sensitivity of the currents $i_d,i_q$ with respect
to the measurement errors. The paper \cite{KSW2021} also presents the
linearization of the model \rfb{eq:ss5}, a typical application example with
relevant Bode plots, and experimental results.

The instantaneous active power $P$ to the grid is \BEQ{eq:P} P \m=\m v_d i_d +
v_q i_q \m=\m -V[i_d\sin\delta+i_q\cos\delta] \end{equation} (see also
\cite[eq.\m(17)]{Natarajan2017}), but this is not computed in the control
algorithm, except possibly for monitoring. It is easy to derive from \rfb{eq:Q}
and \rfb{eq:P} that \BEQ{eq:P_Q} P^2+Q^2 \m=\m V^2(i_d^2 + i_q^2) \m.
\end{equation}

We derive a nice formula linking the $dq$ currents and the powers $P$ and $Q$.
We know from \rfb{eq:Q} and \rfb{eq:P} that \BEQ{Lennon} \bbm{P\\ Q} \m=\m
-V\bbm{ \cos\delta & \sin\delta\\ -\sin\delta & \cos\delta} \bbm{i_q\\ i_d} \m.
\end{equation} By inverting the matrix, we obtain \BEQ{eq:id_iq_P_Q} \bbm{i_q\\
i_d} \m=\m -\frac{1}{V} \bbm{ \cos\delta & -\sin\delta \\ \sin\delta &
\cos\delta} \bbm{P\\ Q} \m. \end{equation}

In Sect. \ref{sec3} we will study the equilibrium points
of the fourth order model \rfb{eq:ss4}, and in Sect. \ref{sec4} we will extend
the study to the fifth order non-saturated model \rfb{eq:ss5}. Finally, the
model \rfb{eq:ss5_sat} will be used in Sect. \ref{sec5} to derive local
exponential stability results for the grid-connected synchronverter.

\section{Equilibrium points of the 4th order grid-connected synchronverter}
\label{sec3} 

In this section we study the equilibrium points of the fourth order model
\rfb{eq:ss4} for the grid-connected synchronverter. Thus, $i_f$ is treated as a
parameter here (i.e., there is no field current controller for the reactive
power $Q$). Our main results is a geometric
representation of the equilibrium points of \rfb{eq:ss4} in the power plane: We
find that, for ``reasonable'' values of $i_f$, the images of the corresponding
equilibrium points of \rfb{eq:ss4} through the mappings $P$ and $Q$ from
\rfb{Lennon} move on a circle in the power plane. The radius and centre of this
circle depend on the synchronverter parameters, on the grid voltage, and on the
prime mover torque at the equilibrium. In addition, the point $(P,Q)$ determines
the power angle $\delta$ at the equilibrium.  Finally, we establish a crucial
results for the stability analysis of Sect. \ref{sec5}: We find the interval of
those field currents $i_f>0$ for which the reactive power $Q$ corresponding to
the relevant equilibrium point is increasing (as a function of $i_f$).

The equilibrium points of \rfb{eq:ss4} have been explicitly computed in
\cite[Sect. 3]{Natarajan2018}, under the assumption of a constant field current
$i_f$ in a reasonable range $I_f\subset(0,\infty)$. For the reader's
convenience, we report those results here. In this paper, angles are regarded
modulo $2\pi$, i.e., $\delta$ and $\delta+2\pi$ are considered to be the same
angle, except for certain arguments in Sect. \ref{sec5}.

\begin{mdframed} \begin{assumption} \label{ass_1} Let
$R,L,J,m,D_p,V,\o_g,\o_n>0$ and $T_m\in\rline$ be given. Denote \begin{equation}
\label{eq:Tm_tilde} \tilde{T}_m \m=\m T_m+D_p(\o_n-\o_g). \end{equation} Assume
that \begin{equation}\label{eq:eq} 4R\o_g\tilde{T}_m\geq-V^2. \end{equation}
\end{assumption} \end{mdframed}

\medskip{\color{blue} \begin{proposition} \label{prop:Xi} Consider the model
\rfb{eq:ss4}, with $\xxx$ from \rfb{eq:x}, and with parameters satisfying
Assumption \ref{ass_1}. Denote \begin{equation} \label{phi} \phi
\in\left(0,\frac{\pi}{2}\right) \ \ \ \text{such that} \ \ \ \tan\phi \m=\m
\frac{\o_gL}{R} \m, \end{equation} \begin{equation} \label{eq:Lambda} \L(i_f)
\m=\m -\frac{\tilde{T}_m}{mi_f} \frac{L\sqrt {p^2+\o_g^2}} {V} + \frac{mi_f\o_g
p}{V\sqrt{p^2+\o_g^2}}, \end{equation} where $p=R/L$. We define the interval
$I_f\subset(0,\infty)$ as follows: \begin{equation*} \label{eq:I} I_f \m=\m
\{i_f > 0\m\big|\m \ |\L(i_f)| \leq 1\}. \end{equation*}

For any $i_f\in I_f\cup(-I_f)$, the model \rfb{eq:ss4} has two equilibrium
points, $\xxx^e_1$ and $\xxx^e_2$, with the power angles $\delta^e_1$ and
$\delta^e_2$ satisfying: \begin{equation} \label{eq:delta_12} \delta^e_1 \m=\m 
\arccos\L-\phi, \qquad \delta^e_2 \m=\m -\arccos\L-\phi, \end{equation} where
$\arccos\L:[-1,1]\to[0,\pi]$. The other components of the equilibrium states
$\xxx_j^e$ are given (for $j\in\{1,2\}$) by \begin{equation} \label{eq:id_iq_w}
i_{dj}^e \m=\m -\frac{\tilde{T}_m\o_g}{m i_f p} + \frac{V\sin \delta^e_j}{R},
\quad i_q^e \m=\m -\frac{\tilde{T}_m}{mi_f}, \quad \o^e \m=\m \o_g.
\end{equation} Note that if $|\L|=1$, then $\delta^e_1=\delta^e_2$ and thus
$\xxx^e_1=\xxx^e_2$. \end{proposition}} \vspace{1mm}

Note that Assumption \ref{ass_1} guarantees that $I_f$ is nonempty. \vspace{1mm}

\begin{figure}[h] \centering 
\includegraphics[width=0.7\linewidth]{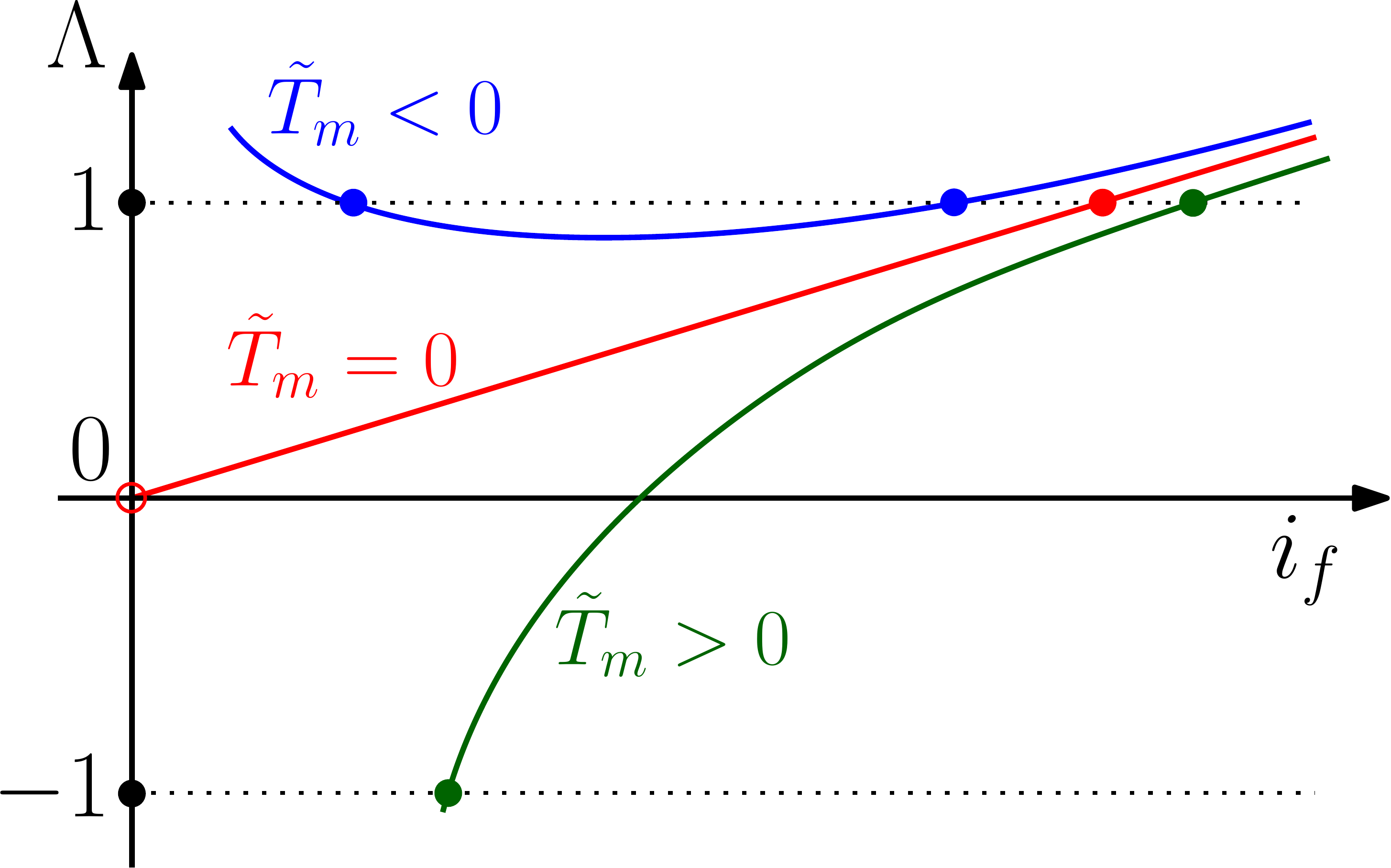}
\caption{The plot of $\L(i_f)$ for different values
of $\tilde{T}_m$ and for $i_f>0$.} \label{fig:Lambda} \end{figure}

It is clear from \rfb{eq:swing} that $\tilde{T}_m$
represents the prime mover torque at equilibrium. The proof follows from
\cite[Sect.~3]{Natarajan2018}, where the notation is slightly different: what is
denoted in \cite{Natarajan2018} by $T_m$, $R_s$ and $L_s$, is denoted here by
$T_m+D_p\o_n$, $R$ and $L$, respectively. Moreover, in \cite{Natarajan2018} it
is assumed that $i_f>0$, however the derivations in
\cite[Sect.~3]{Natarajan2018} remains valid also for $i_f\in(-I_f)$. We now
prove that if $\tilde{T} _m\not=0$, then $I_f$ is a closed interval. If
$\tilde{T}_m>0$, then $\L$ is an increasing function of $i_f$, and our claim
follows. If, instead, $\tilde{T}_m<0$, then $\L$ is first decreasing for a
certain interval of $i_f$, after which it becomes increasing, and we have $\L>0$
for all $i_f>0$. Thus, we can conclude again that $I_f$ is a closed interval.
Finally, if $\tilde{T}_m=0$, then $\L$ depends linearly on $i_f$ and it is clear
that $I_f$ is an interval (not closed). The above scenarios are depicted in
Fig.~\ref{fig:Lambda}.

\begin{remark} \label{FIAT} As mentioned above, if $\tilde{T}_m>0$ then $\L$ is
an increasing function of $i_f>0$ and \vspace{-1mm} $$\{\L(i_f) \m\big|\m i_f\in
I_f\} \m=\m [-1,1].$$ Thus, for every $\L\in[-1,1]$, \rfb{eq:Lambda} has two
solutions: \vspace{-1mm} \begin{equation} \label{eq:i_f1_i_f2} i_{f1} \m=\m
\frac{\sqrt{p^2+\o_g^2}\left(\L V+ \sqrt{\L^2V^2+4\o_g
R\tilde{T}_m}\right)}{2m\o_g p}, \end{equation} and $i_{f2}$ is as above, with
$-$ instead of $+$ in front of the square root in the brackets. Clearly
$i_{f2}<0<i_{f1}$. Thus, there is only one positive solution of \rfb{eq:Lambda}
for each fixed $\L\in[-1,1]$.

On the contrary, if $\tilde{T}_m<0$ then $\L$ is first a decreasing function of
$i_f>0$, and then an increasing one. Moreover, $\L(i_f)>0$ for all $i_f>0$ (see
Fig.~\ref{fig:Lambda}). This implies that \vspace{-1mm} $$\{\L(i_f) \m\big|\m
i_f\in I_f\}\subset(0,1],$$ and that if $\L$ belongs to the above set, then
$i_{f1}, i_{f2}$ from \rfb{eq:i_f1_i_f2} are both positive. Finally, for the
case $\tilde{T}_m=0$, $i_{f2}=0$ and $\L$ is linear in $i_f$, so that $\{\L(i_f)
\m\big|\m i_f\in I_f\}=(0,1]$. \end{remark}

\vspace{1mm}{\color{blue} \begin{proposition} \label{Rouhani}
We use the notation of Proposition \ref{prop:Xi}. Under
Assumption \ref{ass_1}, if for some $i_f\in I_f$ the model \rfb{eq:ss4} has a
stable equilibrium point $\xxx_{st}^e$, then $\xxx_{st}^e=\xxx^e_1$ and
$$\delta^e_1 \m\in\m (-\phi,\pi-\phi) \m.$$ \end{proposition}}

Note that if $R$ tends to zero then $\phi$ tends to $\pi/2$, see \rfb{phi}, and
the above condition becomes the famous necessary stability condition
$\delta\in(-\frac{\pi}{2},\frac{\pi}{2})$ appearing often in the literature.
\vspace{1mm}

{\it Proof.} \m Denote by $h(\xxx)$ the right-hand side of
\rfb{eq:ss4}. Let $A_{lj}$, $j\in\{1,2\}$, be the Jacobian $A_{lj}=\partial
h/\partial\xxx$ computed at $\xxx^e_j$. A necessary condition for the
equilibrium point $\xxx^e_j$ to be stable is that $\HHH^{-1}A_{lj}$ is a stable
matrix, which implies that $\det(\HHH^{-1}A_{lj})>0$. It can be verified easily
that \begin{equation} \label{eq:det_A_l} \det(\HHH^{-1}A_{lj}) \m>\m 0 \iff
\sin(\delta^e_j+\phi)>0 \m, \end{equation} see \cite[eq.\m(3.5)]{Natarajan2018}
for the detailed derivation.

Recall the expressions of $\delta^e_1$ and $\delta^e_2$ from \rfb{eq:delta_12}.
We have $$\delta^e_1+\phi\in[0,\pi], \qquad \delta^e_2+\phi\in[-\pi,0].$$
Clearly, \rfb{eq:det_A_l} holds only for $j=1$. Thus, $\xxx_{st}^e=\xxx^e_1$
and, moreover, we must have $\delta^e_1+\phi\in(0,\pi)$. \hfill $\blacksquare$

\medskip{\color{blue} \begin{proposition} \label{tan_delta}
We use the notation of Proposition \ref{prop:Xi}. Consider the model
\rfb{eq:ss4}, with parameters satisfying Assumption \ref{ass_1}, and let $i_f\in
I_f\cup(-I_f)$, so that \rfb{eq:ss4} has two equilibrium points. Then, at every
equilibrium point of this system we have \begin{equation} \label{eq:Tm_wg4}
\tilde{T}_m \o_g \m=\m P + R \frac{P^2+Q^2}{V^2} \m. \end{equation} Moreover,
the power angle value $\delta^e$ at the equilibrium satisfies \begin{equation}
\label{eq:tan_delta4} \tan \delta^e \m=\m \frac{\o_g L P-RQ} {RP + \o_g LQ +
V^2} \m. \end{equation} 
\end{proposition}}

{\it Proof.} \m In this proof, for convenience, we omit the
superscript $e$ to indicate the equilibrium point values. At an equilibrium
point of \rfb{eq:ss4} the left-hand sides of \rfb{eq:id_new} and \rfb{eq:iq_new}
are zero. We multiply these equations with $i_d$ and $i_q$, respectively, and we
add them using \rfb{eq:v_qv_d}, obtaining $$-R(i_d^2+i_q^2)-mi_f i_q\o - v_d i_d
- v_q i_q \m=\m 0 \m.$$ Using the formulas \rfb{eq:Te} and \rfb{eq:P}, we get
$$R(i_d^2+i_q^2)-T_e\o + P \m=\m 0 \m.$$ It follows from \rfb{eq:swing} that
$T_e=\tilde{T}_m$, and we know from \rfb{eq:id_iq_w} that $\o=\o_g$.
Substituting these values above, we get $$\tilde{T}_m \o_g \m=\m P +
R(i_d^2+i_q^2) \m.$$ Using \rfb{eq:P_Q} this becomes \rfb{eq:Tm_wg4}.

Now we turn our attention to \rfb{eq:tan_delta4}. If we multiply both sides of
\rfb{eq:id_new} (at equilibrium) with $\sin\delta$\m, both sides of
\rfb{eq:iq_new} (at equilibrium) with $\cos\delta$ \m and then we add them, we
get \vspace{-1mm} $$ m\o_gi_f\cos\delta \m= -R[i_d\sin\delta \m + i_q\cos\delta]
\m + \o_g L\frac{Q}{V} \m+ V\m.$$ In the same way, if we multiply
\rfb{eq:id_new} with $\cos \delta$\m, \rfb{eq:iq_new} with $\sin\delta$ and we
subtract them, we get \vspace{-2mm} $$ -m\o_gi_f\sin\delta \m=
\o_gL[i_d\sin\delta \m + i_q\cos\delta] \m + R\frac{Q}{V} \m.$$ According to
\rfb{eq:P} the last two equations can be rewritten as \begin{equation}
\label{eq:eq14} m i_f\o_g\cos\delta \m=\m R\frac{P}{V} + \o_g L\frac{Q}{V} + V
\m, \end{equation} \begin{equation} \label{eq:eq24} m i_f\o_g\sin\delta \m=\m
\o_gL\frac{P}{V} \m - R\frac{Q}{V} \m. \end{equation} Since $i_f\neq0$, the left
sides of \rfb{eq:eq14} and \rfb{eq:eq24} cannot be both zero. We divide the
sides of \rfb{eq:eq24} by the sides of \rfb{eq:eq14}, which shows that
\rfb{eq:tan_delta4} holds. \hfill $\blacksquare$ \vspace{1mm}

\begin{remark} The equation \rfb{eq:Tm_wg4} has a clear intuitive
interpretation: the left-hand side is the mechanical power coming from the
virtual prime mover (the frequency droop mechanism is part of the prime mover).
The second term on the right-hand side is the power consumed in the output
resistors $R$ in series with each of the three phases, if we think of the model
as representing a synchronous machine. (This follows from \rfb{eq:P_Q} and the
fact that the Park transformation is unitary.) \end{remark}

\vspace{1mm} In the following, we present the novel geometric representation of
the (manifold of) equilibrium points of \rfb{eq:ss4} in the $PQ$ plane. For
this, we first introduce some useful notation. \vspace{1mm}

{\bf Notation.} \m We use the notation of Proposition
\ref{prop:Xi}. Consider the model \rfb{eq:ss4}, with parameters satisfying
Assumption \ref{ass_1}. We define $i_{f-}=\inf\m I_f$, and $i_{f+}=\sup\m I_f$.
(Depending on the sign of $\tilde{T}_m$, $\L(i_{f-})$ and
$\L (i_{f+})$ take different values, as discussed in Remark \ref{FIAT}.) Let
$\xxx^e_1(i_f)$, $\xxx^e_2(i_f)$ be the two equilibrium points of \rfb{eq:ss4}
corresponding to $i_f\in I_f$, as described in
\rfb{eq:delta_12}-\rfb{eq:id_iq_w}. ($\xxx^e_1(i_f)$, $\xxx^e_2 (i_f)$ coincide
at $i_f=i_{f-}$ and at $i_f=i_{f+}$.) We denote by
$\delta_j(i_f)$ the power angle component of $\xxx^e_j(i_f)$, $j\in\{1,2\}$, and
by $P_j(i_f)$ ($Q_j(i_f)$) the active (reactive) power at the equilibrium point
$\xxx^e_j(i_f)$, for $j\in\{1,2\}$. If $X,Y,Z\in\rline^2$, then
$\widehat{X;Y;Z}$ denotes the angle from the vector $X-Y$ to the vector $Z-Y$
(counterclockwise). We do not distinguish between a vector and the pair of real
numbers that are its coordinates.

\vspace{1mm}{\color{blue} \begin{theorem} \label{circle} Consider the model
\rfb{eq:ss4}, with parameters satisfying Assumption \ref{ass_1}. Then the points
in $\rline^2$ defined by \vspace{-2mm} $$ S_j(i_f) \m=\m \left(P_j(i_f),\m
Q_j(i_f)\right),\quad i_f\in I_f\m,\ \ j\in\{1,2\}$$ are on the circle with
centre $C$ and radius $r$ given by \begin{equation} \label{eq:r_C} C \m=\m
\left(-\frac{V^2}{2R},\m 0\right)\m ,\qquad r^2 \m=\m
\frac{V^4+4V^2R\tilde{T}_m\o_g}{4R^2} \m. \end{equation} Define the points
$Z,M,O\in\rline^2$ as $$ Z \m=\m \left(R,\m \o_g L\right),\qquad M \m=\m
-\frac{V^2}{\|Z\|^2} Z\m,\qquad O \m=\m (0,0) \m.$$

Then the distances $CO$, $CM$ are equal and \begin{equation} \label{Bennett}
\reallywidehat{O;M;C} \m=\m \reallywidehat{C;O;M} \m=\m \phi \m. \end{equation}
Moreover, the following holds: \begin{equation} \label{eq:M_P_M_Q} S_j(i_f)-M
\m=\m \frac{V}{\|Z\|} \begin{bmatrix} \cos(\phi-\delta_j(i_f)) \\
\sin(\phi-\delta_j(i_f)) \end{bmatrix} m i_f\o_g \m. \end{equation}
\end{theorem}}

{\it Proof.} \m According to Proposition \ref{tan_delta}, the powers $P_j(i_f)$
and $Q_j(i_f)$ satisfy the quadratic equation $$ P_j^2+Q_j^2 + \frac{V^2}{R} P_j
\m=\m \frac{\tilde{T}_m \o_g V^2}{R} \m.$$ The solutions of this equation are on
a circle symmetric with respect to the $P$ axis. The formulas for the centre $C$
and the radius $r$ follow from standard computations.

\begin{figure*}[h!] 
\centering{\subfigure[The case $\o_g L>R$, described in Theorem
\ref{circle_more}\m(a).]
{\includegraphics[height=0.30\textwidth]{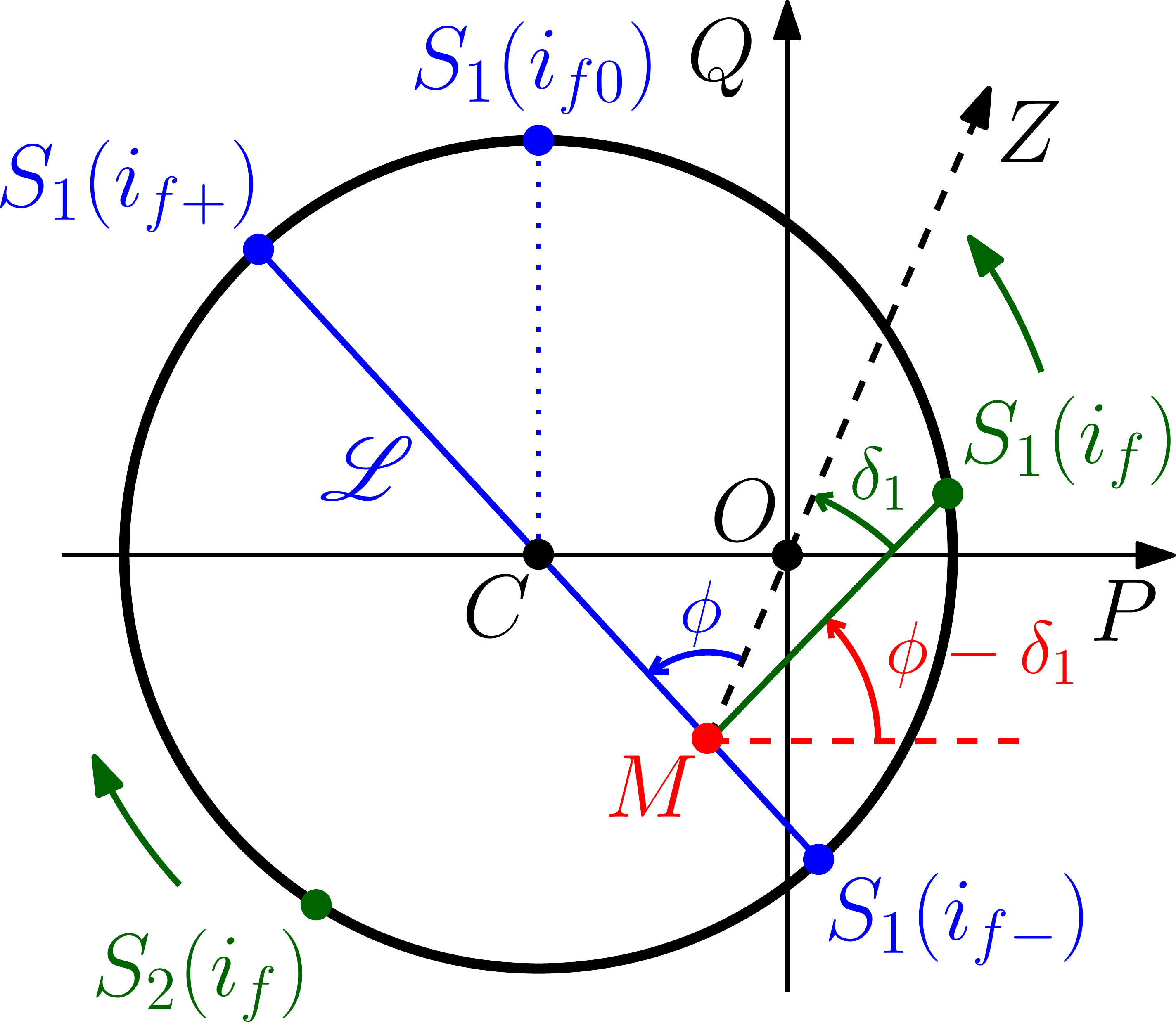}
\label{fig:circle}} \hspace{25mm} \subfigure[The case $\o_g L<R$, described in
Theorem \ref{circle_more}\m(b).]
{\includegraphics[height=0.30\textwidth]{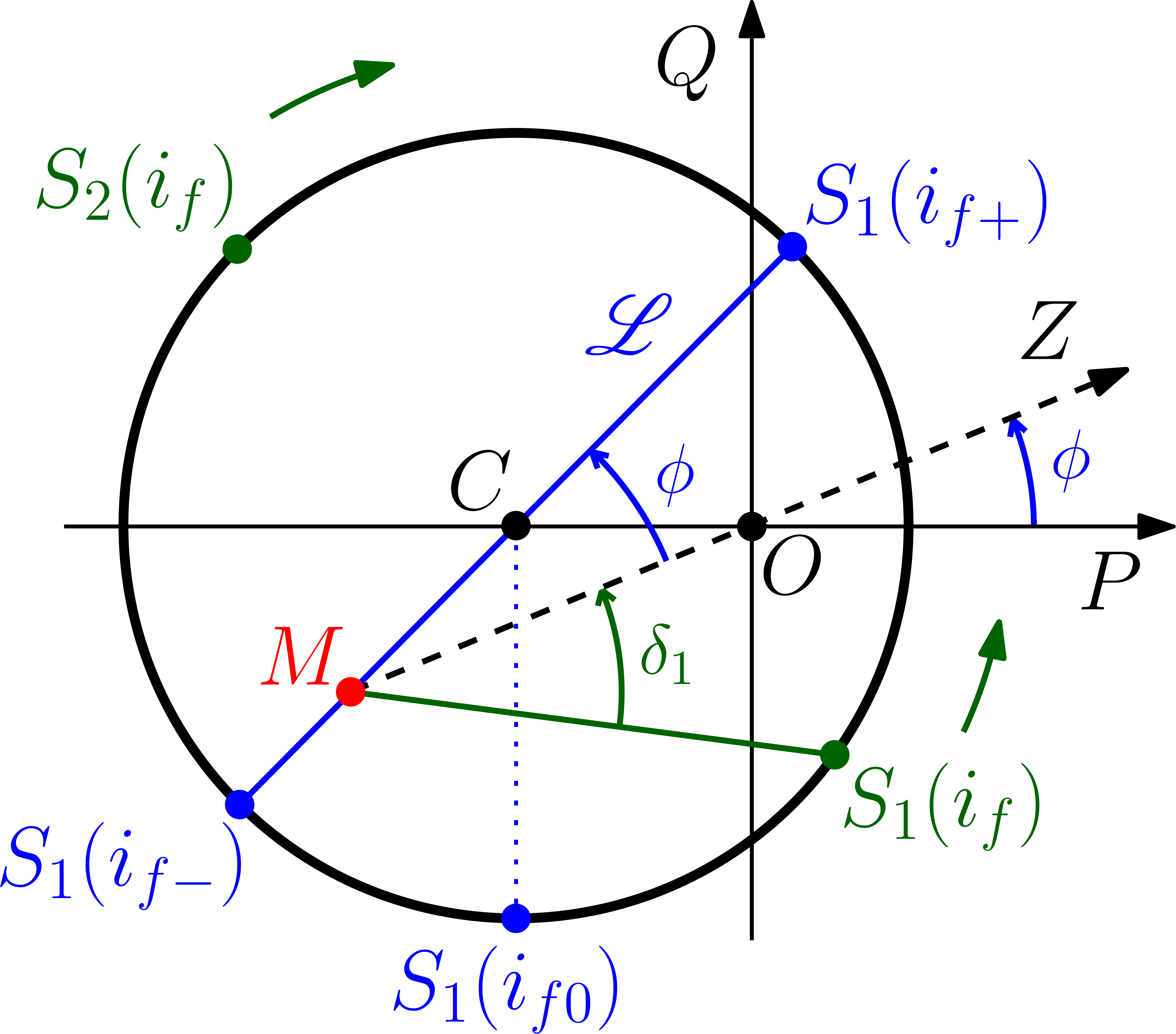} \label{alternative}}}
\caption{The circles (for $\tilde{T}_m>0$) on which the
vectors $S_1(i_f)$ and $S_2(i_f)$ move. The green arrows indicate the movement
of $S_1(i_f)$ and $S_2(i_f)$ for increasing $i_f$.} \end{figure*}

From a routine computation we get that $$\|M-C\| \m=\m \frac{V^2}{2R} \m=\m
\|C\| \m.$$ One conclusion from the above is that the triangle $COM$ is
isosceles, and since the angle of $Z$ (with respect to the $P$ axis) is $\phi$,
we get that \rfb{Bennett} holds (see Fig.~\ref{fig:circle}, \ref{alternative}).

We now prove \rfb{eq:M_P_M_Q}. For convenience we denote
$(P,Q)=(P_j(i_f),Q_j(i_f))$, $\delta=\delta_j(i_f)$. From \rfb{eq:id_new},
\rfb{eq:iq_new} and \rfb{eq:id_iq_P_Q}, we have $$ \frac{1}{V}\bbm{R & -\o L \\
\o L & R}\begin{bmatrix}\cos\delta & \sin\delta \\ -\sin\delta &
\cos\delta\end{bmatrix} \begin{bmatrix} Q \\
P\end{bmatrix}+V\begin{bmatrix}\sin\delta \\
\cos\delta\end{bmatrix}=\begin{bmatrix} 0 \\ mi_f\o_g \end{bmatrix}.$$ Using the
definition of $\phi$ from \rfb{phi}, we have $$ \frac{1}{V}\bbm{R & -\o L \\ \o
L & R}=\frac{\|Z\|}{V}\bbm{\cos \phi & -\sin\phi \\ \sin\phi & \cos\phi}.$$
Substituting this above, commuting the first two matrices,
and multiplying with the inverse of the matrix from \rfb{Lennon}, we obtain $$
\frac{\|Z\|}{V}\bbm{\cos\phi & -\sin\phi \\ \sin\phi & \cos\phi} \begin{bmatrix}
Q \\ P\end{bmatrix}+V\begin{bmatrix}0 \\ 1
\end{bmatrix}=\begin{bmatrix}-\sin\delta \\ \cos\delta \end{bmatrix}mi_f\o_g.$$
Multiplying with the inverse of the first matrix above, and also with $V/\|Z\|$,
and swapping the rows, we get $$ \bbm{P \\ Q}+\frac{V^2}{\|Z\|}\bbm{\cos\phi \\
\sin\phi} = \frac{V} {\|Z\|}\bbm{\cos(\phi-\delta) \\
\sin(\phi-\delta)}mi_f\o_g.$$ From here, we get \rfb{eq:M_P_M_Q} by substituting
$$M=-\frac{V^2}{\|Z\|}\bbm{\cos\phi \\ \sin\phi}.\vspace{-5mm}$$ \hfill
$\blacksquare$

\vspace{1mm} \begin{remark} \label{rmk:star} From \rfb{eq:M_P_M_Q} several
useful facts follow. First, taking the norms, we have that (for $i_f\in I_f$)
\vspace{-1mm} \begin{equation} \label{eq:M_circle}
\left\|S_j(i_f)-M\right\|=\frac{V}{\|Z\|}mi_f\o_g, \end{equation} i.e., the
distance from $S_j(i_f)$ to $M$ is proportional to $i_f$. This implies that the
level curves in the power plane for constant $i_f$ are circles, with centre $M$
and radius given by \rfb{eq:M_circle}. Second, \rfb{eq:M_P_M_Q} tells us that
the vector $S_1(i_f)-M$ forms an angle of $\phi-\delta_1$ with the $P$ axis (see
Fig.~\ref{fig:circle}, \ref{alternative}). Thus, \begin{equation*}
\reallywidehat{S_1(i_f);M;O} \m=\m \phi-(\phi-\delta_1) \m=\m \delta_1 \m.
\end{equation*} Third, clearly $i_{f-}\leq i_f \leq i_{f+}$, for any $i_f\in
I_f$. From \rfb{eq:M_circle}, $S_1(i_{f-})$ is the point on the circle from
Theorem \ref{circle} that is the closest to $M$, while $S_1(i_{f+})$ is the
point on the same circle that is the farthest from $M$. This implies that
$S_1(i_{f-}),\m M,\m C$ and $S_1(i_{f+})$ are on a straight line $\Lscr$, as in
Fig.~\ref{fig:circle}, \ref{alternative}.

From the above facts it follows that, for increasing $i_f\in I_f$, the point
$S_1(i_f)$ moves counterclockwise on the circle described in Theorem
\ref{circle}, from $S_1(i_{f-})$ to $S_1(i_{f+})$. \end{remark}

\vspace{1mm}\begin{remark}\label{rmk:Tm_outside} It follows from the formula for
$r$ in \rfb{eq:r_C} that, depending on the sign of $\tilde{T}_m$, three
scenarios are possible for the points $M$ and $O$ from Theorem \ref{circle}: If
$\tilde{T}_m>0$ then $M,O$ are inside the circle, while if $\tilde{T}_m<0$
($\tilde{T}_m=0$) then $M,O$ are outside (on) the circle. The case
$\tilde{T}_m>0$ is the most common. \end{remark}

\vspace{1mm}{\color{blue} \begin{theorem} \label{circle_more}
We use the notation of Theorem \ref{circle}. Let
Assumption \ref{ass_1} hold. Then: \vspace{1mm}

(a) If $\o_gL>R$, then $M$ is to the right of $C$. There is a unique $i_{f0}\in
I_f$ for which $S_1(i_{f0})=(-V^2/2R,\m r)$ and \vspace{-1mm} $$ \frac{\dd}{\dd
i_f} Q_1(i_f) \m>\m 0 \qquad \mbox{for }\ i_f\m\in\m I_f^+\m=\m(i_{f-},i_{f0})
\m.$$

(b) If $\o_g L\leq R$, then $M$ is to the left of (or directly below) $C$. There
is a unique $i_{f0}\in I_f$ for which $S_1(i_{f0})= (-V^2/2R,\m -r)$ and
\vspace{-2mm} $$ \frac{\dd}{\dd i_f} Q_1(i_f) \m>\m 0 \qquad \mbox{for }\
i_f\m\in\m I_f^+\m=\m(i_{f0},i_{f+}) \m.$$ \end{theorem}}

\vspace{1mm} Theorems \ref{circle}, \ref{circle_more} are illustrated in
Fig.~\ref{fig:circle}, \ref{alternative}. As discussed in
Remark \ref{rmk:star}, we see that $S_1(i_f)$ moves counterclockwise on the
circle, for increasing $i_f$, from $S_1(i_{f-})$ to $S_1(i_{f+})$, while
$S_2(i_f)$ moves clockwise between the same two endpoints. The movement of
$S_2(i_f)$ is symmetric to the one of $S_1(i_f)$, with respect to the line
$\Lscr$.

Note that the case $\o_g L>R$ is the most common.

\medskip {\it Proof of Theorem} \ref{circle_more}. Let $\tilde{T}_m\geq0$. In
the case (a), an elementary computation shows that the $P$-coordinate of $M$ is
larger than that of $C$: \vspace{-1.5mm} $$ -\frac{V^2 R}{\|Z\|^2} \m>\m
-\frac{V^2}{2R} \m.$$ Hence, $M$ is to the right of $C$, as stated. Note that
this implies that the slope of $\Lscr$ is negative, as in Fig. \ref{fig:circle}.

As discussed in Remark \ref{rmk:star}, $S_1(i_f)$ moves counterclockwise on the
circle for increasing $i_f\in I_f$. Since $i_f$ is proportional to the distance
from $S_1$ to $M$, $Q_1$ is strictly increasing (with positive derivative) for
$i_f\m\in\m I_f^+ =(i_{f-},i_{f0})$, where $i_{f0}$ is the field current for
which $Q_1(i_f)$ reaches its maximum value, namely $r$. From Fig.
\ref{fig:circle} we see that $i_{f0}$ is the unique field current for which
$P_1(i_{f0})=-V^2/2R$.

We move now to case (b). We perform the same elementary computation as before,
reaching the opposite conclusion for $\o_g L\leq R$, namely, that $M$ is to the
left of (or directly below) $C$. Thus, for $\o_g L<R$, the slope of $\Lscr$ is
positive (as depicted in Fig. \ref{alternative}), and for $\o_g L=R$, $\Lscr$ is
vertical.

In the proof of (b), the interval on which $Q_1$ is increasing is from $i_{f0}$,
where $Q_1$ is at its minimum, until $i_{f+}$. We see from
Fig.~\ref{alternative} that $S_1(i_{f0})=(-V^2/2R,-r)$.

The proof of the case $\tilde{T}_m<0$ is similar. \hfill $\blacksquare$

%
%

\begin{figure}[h] \centering 
\vspace{-1mm} \includegraphics[width=0.5\linewidth]{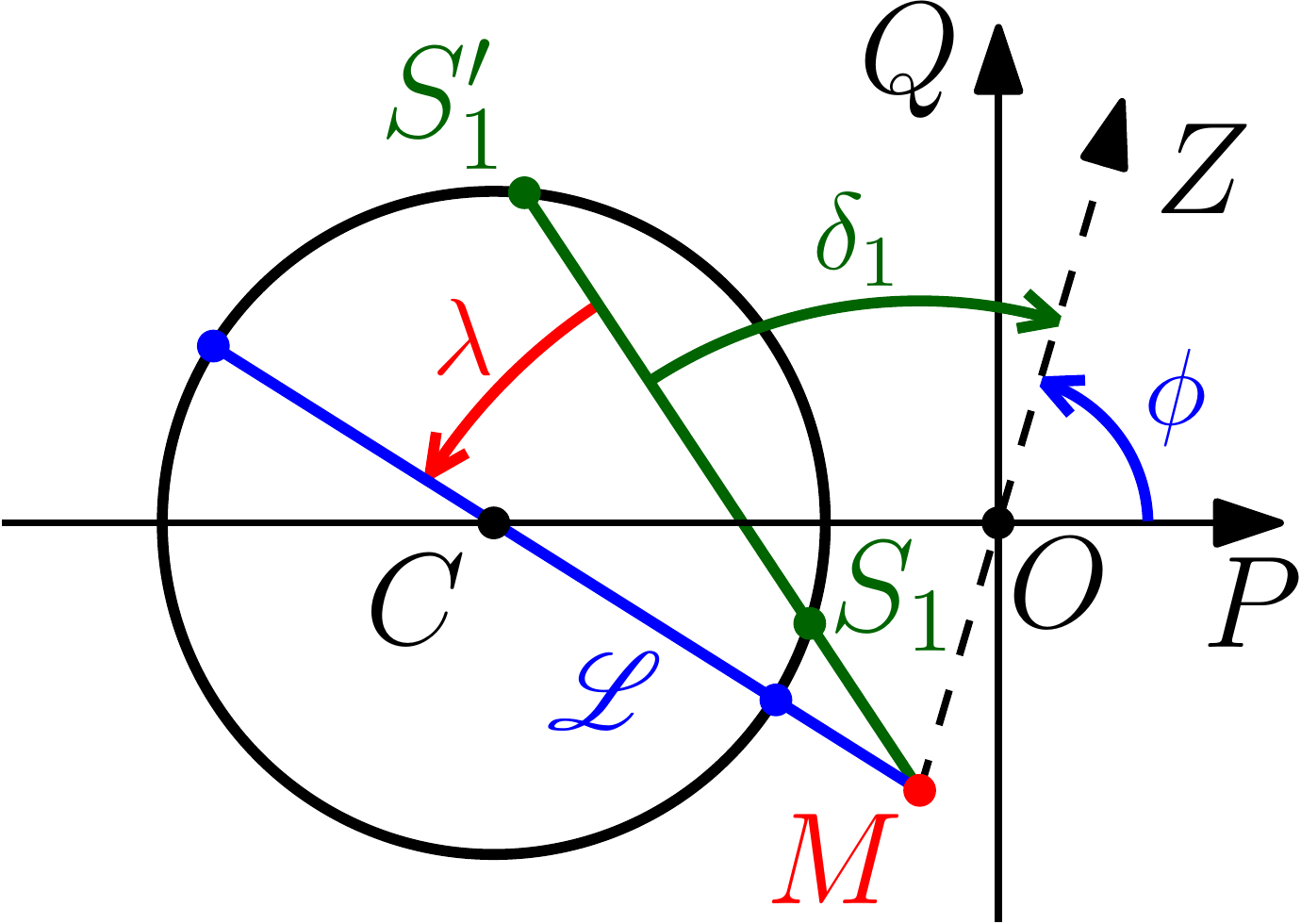}
\vspace{-1mm} \caption{The circle corresponding to
$\tilde{T}_m<0$ and $\o_gL>R$. We have denoted $S_1^{\m \prime}=S_1(i_{f1})$ and
$S_1=S_1 (i_{f2})$ ($0<i_{f2}\leq i_{f1}$). The case $\o_g L\leq R$ is similar,
but derived according to Fig.~\ref{alternative}.} \label{fig:circle_sm}
\end{figure}

\vspace{1mm} \begin{remark} \label{rmk:Tm<0_1} For $\tilde{T}_m<0$ both the
solutions $i_{f1},i_{f2}$ of \rfb{eq:Lambda}, with $\Lambda\in\{\L(i_f)
\m\big|\m i_f\in I_f\}$, are positive (see Remark \ref{FIAT}). This has an
intuitive geometrical meaning. Fixing $\L$ is similar to fixing $\lambda=\arccos
\L$, i.e., the angle $\reallywidehat{S_1(i_{f1});M;C}$ in
Fig.~\ref{fig:circle_sm}. Since $M$ is outside of the circle (see Remark
\ref{rmk:Tm_outside}), the line passing through the points $M$ and $S_1(i_{f1})$
cuts the circle in another point, namely, $S_1(i_{f2})$. The values
$i_{f1},i_{f2}$ are the two (positive) solutions of \rfb{eq:Lambda} mentioned
above. The case $S_1\equiv S_1^{\m \prime}$ in Fig.~\ref{fig:circle_sm}
corresponds to the value of $\L$ for which the square root in \rfb{eq:i_f1_i_f2}
is zero, i.e., $i_{f1}=i_{f2}$. \end{remark}

\vspace{1mm} \begin{remark} \label{rmk:Tm<0_2} The point $S_1(i_f)$ moves
counterclockwise on the circle described in Theorem \ref{circle}, from
$S_1(i_{f-})$ to $S_1(i_{f+})$ for increasing $i_f\in I_f$, as discussed in
Remark \ref{rmk:star}. When $\tilde{T}_m\geq0$, this implies that
$\delta_1(i_f)$ is decreasing from $\delta_1(i_{f-})=\pi-\phi$, to
$\delta_1(i_{f+})= -\phi$. However, when $\tilde{T}_m<0$ then this is not true.
Indeed $\delta_1(i_{f-})=\delta_1(i_{f+})=-\phi$. \end{remark}

\section{Equilibrium points of the fifth order grid-connected synchronverter}
\label{sec4} 

In this section we study the equilibrium points of the
fifth order grid-connected synchronverter model \rfb{eq:ss5}. Using the results
for the fourth order model \rfb{eq:ss4} from the previous section, we derive a
necessary and sufficient condition for the existence of the equilibrium points
of \rfb{eq:ss5} (where $i_f$ is a state variable) and we compute them
explicitly. As in Sect. \ref{sec3}, we consider the grid to be an infinite bus,
with constant $V,\m\o_g$.

\begin{figure}[h] \centering 
\vspace{-2mm} \includegraphics[width=0.9\linewidth]{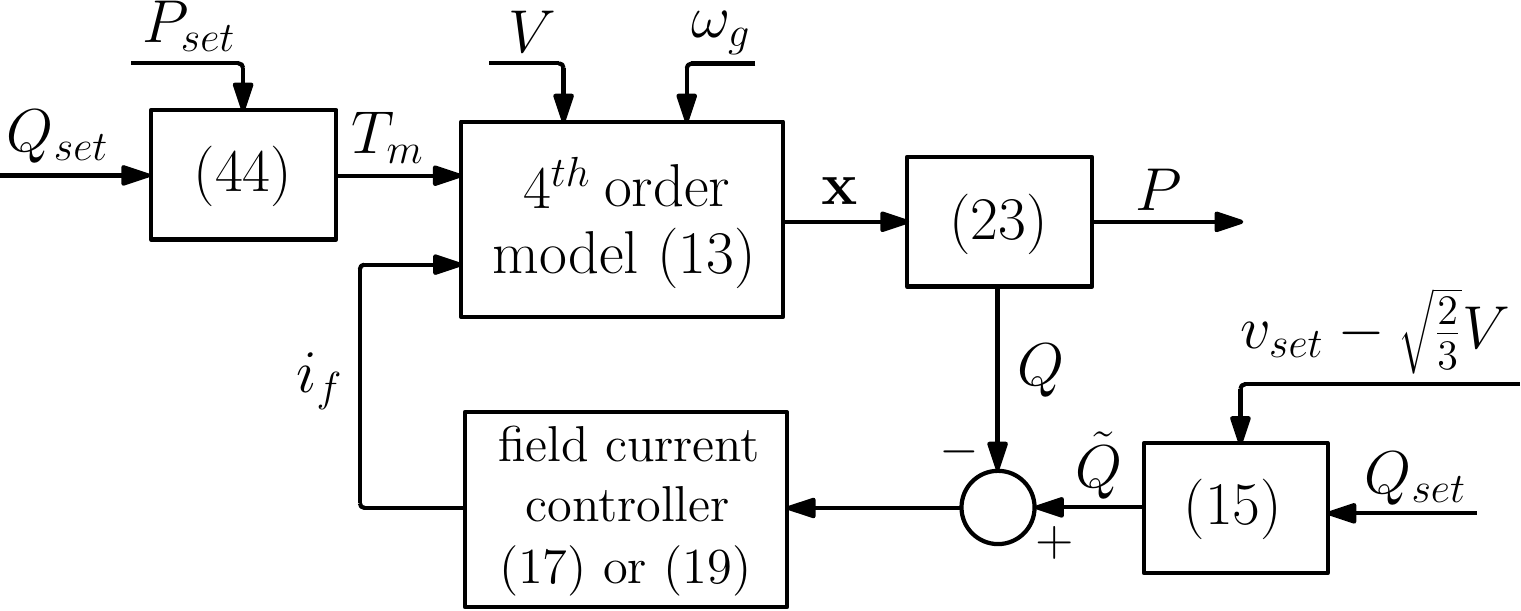}
\vspace{-2mm} \caption{The fourth order model
\rfb{eq:ss4} as a subsystem of the fifth order models \rfb{eq:ss5_sat},
\rfb{eq:ss5}. Depending on the usage of \rfb{eq:Mfif_sat} or \rfb{eq:mif} in
the field current controller, we get respectively \rfb{eq:ss5_sat} or
\rfb{eq:ss5}.} \label{fig:relation} \vspace{-2mm} \end{figure}

The fifth order model \rfb{eq:ss5_sat} or \rfb{eq:ss5} is shown as a block
diagram in Fig.~\ref{fig:relation}, with the fourth order model \rfb{eq:ss4} as
a block.

\begin{mdframed} \begin{assumption} \label{ass_2} Let
$R,L,J,m,D_p,D_q,V,\o_g,\o_n,v_{\rm set}>0$ and $T_m,Q_{\rm set}\in\rline$ be
given. \end{assumption} \end{mdframed}

Our first result concerns mainly the equation that must be satisfied by the
active power $P$ at an equilibrium point of \rfb{eq:ss5}.

\vspace{1mm}{\color{blue} \begin{proposition} \label{prop:3.1} Consider the
model \rfb{eq:ss5}, with parameters satisfying Assumption \ref{ass_2}. Recall
$\tilde{Q}$ from \rfb{eq:Q_tilde} and $\tilde{T}_m$ from \rfb{eq:Tm_tilde}.

A necessary condition for this system to have equilibrium points is
\vspace{-1mm} \BEQ{eq:E} 4R^2\tilde Q^2 \m\leq\m V^4 + 4 R V^2\tilde{T}_m\o_g
\m. \end{equation} At every equilibrium point of this system we have
\BEQ{eq:w_Te_Q} \o \m=\m \o_g \m,\qquad T_e \m=\m \tilde{T}_m \m,\qquad Q \m=\m
\tilde Q, \end{equation} and $P$ satisfies the equation \BEQ{eq:Tm_wg}
\tilde{T}_m \o_g \m=\m P + R \frac{P^2+\tilde Q^2}{V^2} \m. \end{equation}
\end{proposition}}

\begin{remark} \label{Biden} A formula equivalent to
\rfb{eq:Tm_wg} has appeared in \cite[eq.\m(24)]{Natarajan2017}, but instead of a
mathematical proof it was derived from a physical balance equation. As proposed
in \cite{Natarajan2017}, this formula can be used in the synchronverter
algorithm to determine the value of the parameter $T_m$, if the reference values
$P_{\rm set}$ and $Q_{\rm set}$ are given and if some estimate (for instance,
zero) is adopted for the differences $\o_n-\o_g$ and $v_{\rm set}-\sqrt{2/3}V$.
Indeed, if the estimate zero is adopted for these differences (which is, a
priori, our best guess), then \vspace{-1mm} \begin{equation} \label{eq:T_m} T_m
\o_n \m=\m P_{\rm set}+ R \frac{P_{\rm set}^2+Q_{\rm set}^2} {V^2} \m.
\end{equation} For this reason, it is similar to assume that $T_m,\m Q_{\rm
set}$ are given (as in Proposition \ref{prop:3.1}) or that $P_{\rm set},\m
Q_{\rm set}$ are given. \end{remark}

\vspace{1mm} Note that \rfb{eq:E} is equivalent to $|\tilde{Q}|\leq r$, where
$r$ is the radius of the circle from Proposition \ref{circle}. Indeed,
$|\tilde{Q}|>r$ would be an infeasible requirement, as is clear from
Fig.~\ref{fig:circle}, \ref{alternative}.

\vspace{1mm}{\it Proof.} \m We omit the superscript $e$ to
indicate the equilibrium point values. If the system is at an equilibrium point,
then from \rfb{eq:delta} we see immediately that $\o=\o_g$, from \rfb{eq:swing}
we see that \m $T_e=\tilde{T}_m$ and from \rfb{eq:Mfif} we see that $Q=\tilde
Q$\m. Thus, we have proved all the parts of \rfb{eq:w_Te_Q}.

Equation \rfb{eq:Tm_wg} follows from \rfb{eq:Tm_wg4}, substituting $Q=\tilde{Q}$
from \rfb{eq:w_Te_Q}. Note that \rfb{eq:Tm_wg} is a second order equation in
$P$, where the coefficients depend on the parameters of the system. For this
equation to have a real solution, by elementary algebra, the condition
\rfb{eq:E} must be satisfied. Hence, if \rfb{eq:E} does not hold, then the
system cannot have equilibrium points. \hfill $\blacksquare$

\vspace{1mm}\begin{remark} \label{Dr_Fauci} The equilibrium points of
\rfb{eq:ss5} come in symmetric pairs. Indeed, if $\zzz^e=[i^e_d \ i^e_q \ \o_g \
\delta^e \ i^e_f]^\top$ is such an equilibrium point, then also $$
\tilde{\zzz}^e \m=\m [-i^e_d \ -i^e_q \ \o_g \ \delta^e+\pi \ -i^e_f]^\top$$ is
an equilibrium point. The intuition behind this is clear: if we rotate the rotor
by $180^\circ$ and at the same time invert the current $i_f$ in the rotor, then
by the symmetry of the rotor nothing has really changed. The replacement of the
rotor angle $\theta$ with $\theta+\pi$ causes $i_d$ and $i_q$ to change sign,
while the currents in the stationary frame remain unchanged. We see from
\rfb{Lennon} that the active and reactive powers $P,Q$ at $\zzz^e$ and at
$\tilde{\zzz}^e$ are the same. \end{remark}

\vspace{1mm}\begin{remark} \label{Angela} The system \rfb{eq:ss5} has an
exceptional set of equilibrium points corresponding to the point $M$ defined in
Theorem \ref{circle}. Indeed, when the circle defined in Theorem \ref{circle}
passes through the point $M$ (this happens for $\tilde{T}_m=0$), and the values
of $P$ and $\tilde{Q}$ are the coordinates of $M$, namely \BEQ{Boris_Johnson} P
\m=\m -\frac{V^2 R}{R^2+\o_g^2 L^2} \m,\qquad \tilde Q \m=\m -\frac{V^2 \o_g
L}{R^2+\o_g^2 L^2} \m, \end{equation} then if we choose $i^e_f=0$ and any angle
$\delta^e$, we get an equilibrium point of \rfb{eq:ss5}.
This can be checked through a somewhat tedious computation (using
\rfb{eq:id_iq_P_Q}), which shows that for $i^e_f=0$ and any $\delta^e$,
\rfb{eq:id_new} and \rfb{eq:iq_new} hold with zero on the left-hand side. The
other equilibrium equations are easily seen to hold. Thus, for $\tilde{T}_m=0$
and $P, \m\tilde Q$ as in \rfb{Boris_Johnson} we have infinitely many
equilibrium points.

The physical interpretation of these equilibrium points is as follows: here the
rotor has no current and hence no magnetic field, so that its angle is
irrelevant for what happens in the stator windings. The SG now consists of only
the stator windings connected to the power grid, consuming power.
The practical importance of the exceptional set of
equilibrium points discussed above is very small, along with all the equilibrium
points that correspond to negative $i_f$. Indeed, the actual field current
controller employs a saturating integrator (see \rfb{eq:Mfif_sat}), which
constrains the $i_f$ values to an interval of positive numbers (contained in
$I_f$). This is a safety feature that prevents the system from leaving its
normal operating range. \end{remark}

\vspace{1mm} {\color{blue} \begin{theorem} \label{Fakhrizadeh}
We work with the notation of Proposition \ref{prop:3.1}.
Then the model \rfb{eq:ss5}, with parameters satisfying Assumption \ref{ass_2},
has equilibrium points if and only if \rfb{eq:E} is satisfied. Suppose that the
condition \rfb{eq:E} is true, and let us denote by $P_l$ and $P_r$ the two real
solutions of \rfb{eq:Tm_wg}, so that $P_l\leq P_r$, and
$\frac{P_l+P_r}{2}=-\frac{V^2}{2R}$. At every equilibrium point we have $P=P_l$
or $P=P_r$.

Recall the exceptional point $M$ discussed in the last remark. Assume that the
equilibrium point is such that $(P,\tilde Q)\not=M$. Then the angle $\delta^e$
satisfies \vspace{-2mm} \BEQ{eq:tan_delta} \tan \delta^e \m=\m \frac{\o_g L
P-R\tilde Q}{RP + \o_g L \tilde Q + V^2} \m. \end{equation}

If the angle $\delta$ is measured modulo $2\pi$, and \rfb{eq:E} holds with
strict inequality, then the model \rfb{eq:ss5} has precisely four equilibrium
points. Two of them, denoted by $\zzz^e_r$ and $\zzz^e_l$, have the property
that $i_f^e>0$. At $\zzz^e_r$, $P=P_r$, and at $\zzz^e_l$, $P=P_l$. There are
also the two symmetric equilibrium points $\tilde{\zzz}_r^e$ and
$\tilde{\zzz}_l^e$ where $i_f^e<0$, as described in Remark \ref{Dr_Fauci}. If
\rfb{eq:E} holds with equality, then $P_l=P_r=-V^2/2R$ and the model has
precisely two equilibrium points, which are a symmetric pair. \end{theorem}}

\vspace{1mm} \begin{remark} We see from \rfb{eq:tan_delta}
that to any $(P,\tilde Q)\not=M$ in the power plane correspond two possible
equilibrium angles, that differ by $\pi$. This is true also if the denominator
is zero, in that case $\delta^e=\pm\pi /2$. For the exceptional pair $M$, the
right-hand side of \rfb{eq:tan_delta} is $0/0$, so that $\delta^e$ could take
any value, in accordance with Remark \ref{Angela}. \end{remark}

\vspace{1mm} {\em Proof.} \m We omit the superscript $e$ to indicate the
equilibrium point values. Assume that \rfb{eq:E} holds, so that \rfb{eq:Tm_wg}
has two real solutions, $P_l$ and $P_r$, with $P_l\leq P_r$. We know from
Proposition \ref{prop:3.1} that at every equilibrium point, $P=P_l$ or $P=P_r$.

Equation \rfb{eq:tan_delta} follows from \rfb{eq:tan_delta4}, substituting
$Q=\tilde{Q}$ from \rfb{eq:w_Te_Q}. For each choice of $P$ (either $P_l$ or
$P_r$) such that $(P,\tilde Q)\not=M$, this equation has precisely two solutions
modulo $2\pi$, that differ by an angle of $\pi$. (In the extreme case when the
denominator in \rfb{eq:tan_delta} is zero, then the solutions are $\pm\pi/2$.)

Suppose that \rfb{eq:E} holds with strict inequality, which implies that
$P_l<P_r$, and suppose that $(P,\tilde Q)\not=M$. Then we obtain four candidate
equilibrium angles $\delta$ (two for $P=P_l$ and two for $P=P_r$). We now show
that each of these four angles actually corresponds to an equilibrium point
$\zzz=(i_d,i_q,\o_g,\delta,i_f)$. From \rfb{eq:id_iq_P_Q} we see that at any
equilibrium point $$ \bbm{i_q\\ i_d} \m=\m -\frac{1}{V} \bbm{ \cos\delta & -\sin
\delta\\ \sin\delta & \cos\delta} \bbm{P\\ \tilde Q} \m,$$ where $P=P_l$ or
$P=P_r$. From \rfb{eq:Te} and \rfb{eq:w_Te_Q} we see that at any equilibrium
point, \vspace{-2mm} \begin{equation*} \tilde{T}_m \m=\m -m i_f i_q \m.
\end{equation*} Thus, if $\tilde{T}_m\not=0$, $i_f$ can be computed from here.
If $\tilde{T}_m=0$, then \rfb{eq:iq_new} (at the equilibrium) should be used
instead, as long as $(P,\tilde Q)\not=M$. The exceptional case when $(P,\tilde
Q)=M$ leads to $i_f=0$ and arbitrary $\delta$, as discussed in Remark
\ref{Angela}.

It is easy to see that the points $\zzz=(i_d,i_q,\o_g,\delta,i_f)$ computed as
described are indeed equilibrium points, and they come in two symmetric pairs,
as described in Remark \ref{Dr_Fauci}.

When we have equality in \rfb{eq:E}, then $P_l=P_r=-V^2/2R$. Correspondingly,
there are only two solutions for \rfb{eq:tan_delta} (modulo $2\pi$) and they
differ by $\pi$. The currents $i_d,\m i_q,\m i_f$ are computed as before, and we
obtain two equilibrium points (a symmetric pair), one with $i_f>0$ and the other
one with $i_f<0$. \hfill $\blacksquare$

\vspace{1mm} \begin{remark} Under the conditions of the
last theorem, it is easy to see that $P_r\geq 0$ if and only if \BEQ{eq:RQ} R
\tilde Q^2 \m\leq\m V^2 \tilde{T}_m \o_g \m, \end{equation} and $P_r=0$ if and
only if we have equality in \rfb{eq:RQ}. Note that \rfb{eq:RQ} implies
\rfb{eq:E} and we always have $P_l<0$. If an equilibrium point corresponds to
$P_r>0$ and $\tilde Q=0$, then $\tan\delta^e>0$ (this means that
$\delta^e\in(0,\pi/2)\cup(\pi, 3\pi/2)$). Indeed, this can be seen directly from
\rfb{eq:tan_delta}. (These facts are clear from Fig.~\ref{fig:circle},
\ref{alternative}.) \end{remark}

\vspace{1mm} \begin{remark} As mentioned at the end of Remark \ref{Angela}, the
real system \rfb{eq:ss5_sat} can never reach the two equilibrium points with
$i_f^e\leq0$, due to the saturating integrator used in the field current
controller (see \rfb{eq:Mfif_sat}). \end{remark}

\section{Stability of the grid-connected synchronverter} \label{sec5} 

In this section we investigate the stability of the grid-connected
synchronverter model \rfb{eq:ss5_sat} using \cite[Theorem 4.3] {Lorenzetti2020},
which is based on singular perturbation theory. Our main
result in Theorem \ref{thm:5.1} proves that, under reasonable assumptions, there
exists a $\kappa>0$ such that if $\tilde{K}> \frac{1}{\kappa}$, then the fifth
order model \rfb{eq:ss5_sat} has a (locally) exponentially stable equilibrium
point with a ``large'' domain of attraction. This stable equilibrium point
``corresponds'' to $\xxx_1^e$ from Proposition \ref{prop:Xi}. After stating our
main result, we offer a visual representation of the stability region of the
fifth order model \rfb{eq:ss5_sat}, based on the geometric description
introduced in Proposition \ref{circle}.

Note that results closely related to our Theorem \ref{thm:5.1}, with most of the
proof missing, assuming that the model \rfb{eq:ss4} is almost globally
asymptotically stable for every constant $i_f\in [u_{min},u_{max}]$, have been
presented in \cite[Theorem 5.1]{Natarajan2017}. 

We introduce a function $\Xi$ that maps ``reasonable'' values of $i_f$ into the
corresponding first equilibrium point $\xxx^e_1$ of the fourth order model
\rfb{eq:ss4} (see Proposition \ref{prop:Xi}) as \vspace{-1mm} \begin{equation}
\label{eq:Xi} \Xi:I_f\to\rline^4 \quad \text{such that} \quad \Xi(i_f) \m=\m [\m
i_{d1}^e \ i_q^e \ \o_g \ \delta^e_1\m]^\top, \end{equation} where $i_{d1}^e$,
$i_q^e$ and $\delta^e_1$ are given by \rfb{eq:delta_12}-\rfb{eq:id_iq_w}, so
that $\xxx^e_1=\Xi(i_f)$. Here angles are not identified modulo $2\pi$, because
we use results from singular perturbations theory that have been formulated for
systems evolving on $\mathbb{R}^n$. We consider $\delta_1^e\in[-\pi,\pi]$.

Recall the interval $I_f^+$ from Theorem \ref{circle_more}. Then it follows from
the just mentioned theorems that if \rfb{eq:eq} holds with strict inequality (so
that $I_f^+$ is nonempty), then\vspace{-2mm} $$ \frac{\dd}{\dd i_f} Q_1(i_f)
\m>\m 0 \qquad \mbox{for }\ i_f\in I_f^+ \m.\vspace{-2mm}$$

Let $\zzz^e_r=[\m i^e_{dr} \ i^e_{qr} \ \o_g \ \delta^e_r \ i^e_{fr}]^\top$ be
defined as in Theorem \ref{Fakhrizadeh} (i.e., $\zzz^e_r$ is the equilibrium
point of the fifth order model \rfb{eq:ss5} at which $i^e_{fr}>0$ and $P=P_r$).
Assume that $i_{fr}^e \in I_f^+$\m. We see from Fig.~\ref{fig:circle} and Fig.
\ref{alternative} that this implies that the point $(P_1(i^e_{fr}),
Q_1(i^e_{fr}))$ is to the right of the line $\Lscr$ in the power plane, so that
$P_1(i^e_{fr})=P_r$ and thus \vspace{-1mm} $$ \zzz^e_r \m=\m \bbm{\Xi(i_{fr}^e)
\\ i_{fr}^e}.$$

\begin{figure}[t] \centering 
\includegraphics[width=0.6\linewidth]{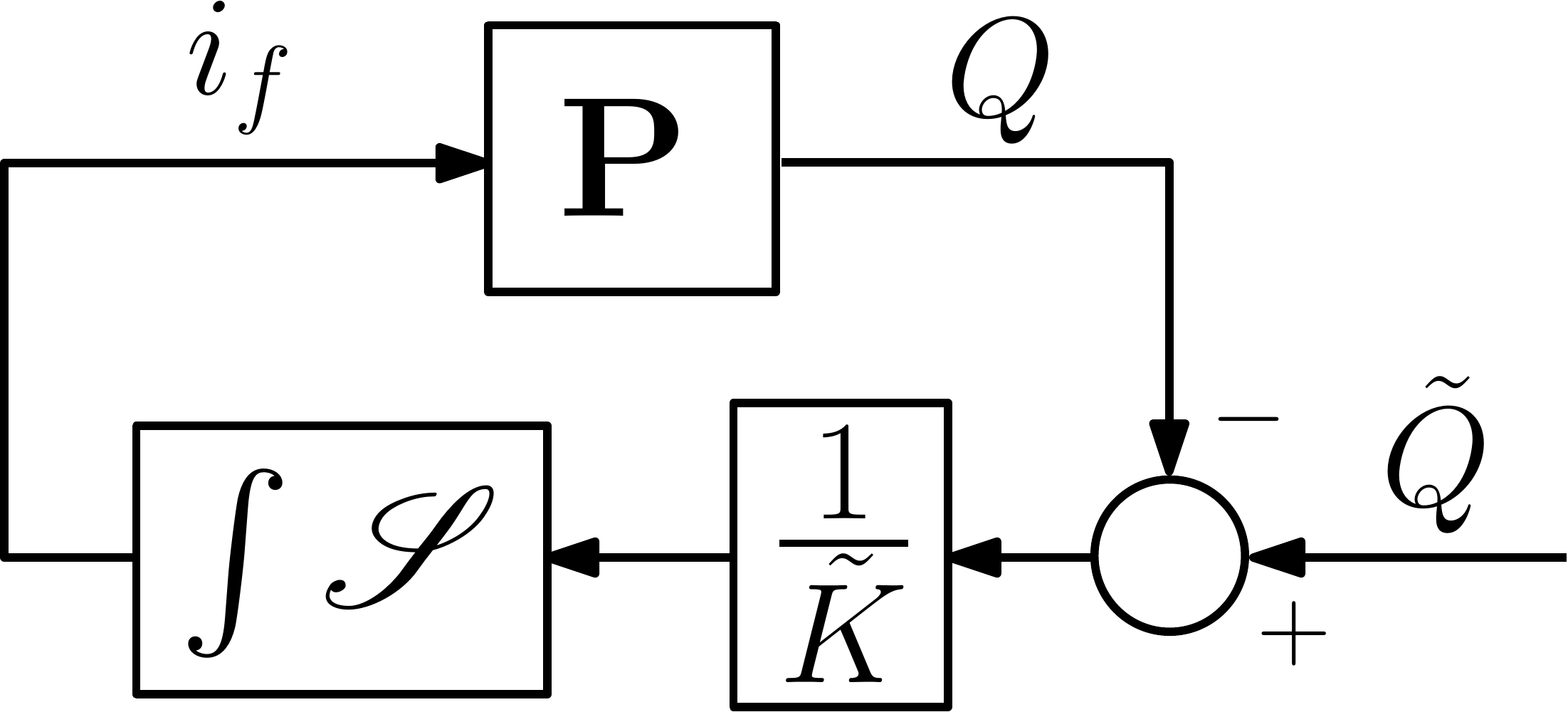} \vspace{-1mm}
\caption{The closed-loop system \rfb{eq:ss5_sat} formed by the plant $\PPP$
from \rfb{eq:ss4}, $\tilde{Q}$ from \rfb{eq:Q_tilde}, and the saturating
integrator from \rfb{eq:Mfif_sat}.} \label{fig:cl} \end{figure}

\vspace{1mm} {\color{blue} \begin{proposition} \label{prop:4.1} We consider the
fourth order system $\PPP$ described by \rfb{eq:ss4}, with parameters satisfying
Assumption \ref{ass_1} (\rfb{eq:eq} with strict inequality). Recall the function
$\Xi$ from \rfb{eq:Xi}. Let $u_{min}<u_{max}$ in $\rline$ and $\e>0$ be such
that $U_\e=[u_{min}-\e,u_{max}+\e]\subset I_f^+$. Denote by $G:U_\e\to\rline$
the steady-state input-output map associated to $\PPP$, with input $i_f$ and
output $Q_1(i_f)$, i.e., $G(i_f)$ is the output $Q_1(i_f)$ at
$\xxx=\Xi(i_f)=\xxx_1^e$. Then, \vspace{-1mm} $$G'(i_f)>0, \FORALL i_f\in
U_\e.\vspace{-2mm}$$ \end{proposition}}


\vspace{1mm} {\em Proof.} \m It follows from Theorem \ref{circle_more} that
$Q_1(i_f)$ is increasing for $i_f\in I_f^+$. Thus, $G'(i_f)>0$ for all $i_f\in
U_\e\subset I_f^+$. \hfill $\blacksquare$

\vspace{1mm} {\color{blue} \begin{theorem} \label{thm:5.1} Consider the model
\rfb{eq:ss5_sat}, with given $R,L,J,m,D_p,D_q,V, \o_g,\o_n,v_{\rm set}>0$ and
$T_m\in\mathbb{R}$, and with the state $\zzz=\sbm{{\bf x}\\ i_f}\in\rline^5$
($\xxx$ is as in \rfb{eq:x}). We use the notation
$\PPP,\m\tilde{T}_m,\m I_f,\m I_f^+,\m u_{min},\m u_{max},\m \e,\m U_\e,\m
\Xi,\m G$ as in Proposition \ref{prop:4.1}. Assume that \rfb{eq:eq} holds with
strict inequality and that the synchronverter parameters are chosen so that
$\PPP$ has a locally exponentially stable equilibrium point for every $i_f\in
U_\e$.

Then, for any $\tilde{Q}\in[G(u_{min}),G(u_{max})]$, denoting
$i_{fr}^e=G^{-1}(\tilde{Q})$, there exist an $\e_0>0$ and a $\kappa>0$ such
that: If $\tilde{K}>\frac{1}{\kappa}$, then $\zzz_r^e= (\Xi(i_{fr}^e),i_{fr}^e)$
is a (locally) exponentially stable equilibrium point of the closed-loop system
\rfb{eq:ss5_sat}, with state space $X=\rline^4\times[u_{min},u_{max}]$.
Moreover, if the initial state $(\xxx(0),i_f(0))\in X$ of \rfb{eq:ss5_sat}
satisfies $\|x(0)-\Xi(i_f(0))\|\leq\e_0$, then $$ \xxx(t)\rarrow\Xi(i_{fr}^e),
\qquad i_f(t)\rarrow i_{fr}^e, \qquad Q(t) \rarrow \tilde{Q},$$ and this
convergence is at an exponential rate. \end{theorem}}

\vspace{1mm} {\em Proof.} The exponential stability of $\PPP$ for each $i_f\in
U_\e$ (as assumed in the theorem) implies the uniform exponential stability of
$\PPP$, see \cite[Remark 3.1]{Lorenzetti2020}. This, together with the result
from Proposition~\ref{prop:4.1}, allows us to apply
\cite[Theorem~4.3]{Lorenzetti2020}, completing the proof of the theorem. \hfill
$\blacksquare$

\vspace{1mm} \begin{remark} The local exponential stability assumption in the
above theorem is true if the parameters satisfy the numerical conditions
presented in \cite[Theorem 1]{Barabanov2017} or in \cite[Theorem
6.3]{Natarajan2018} (the conditions in these two references are not equivalent).
Actually, \cite{Barabanov2017} and \cite{Natarajan2018} conclude aGAS.
\end{remark}

\begin{figure}[h] \centering 
\includegraphics[width=0.8\linewidth]{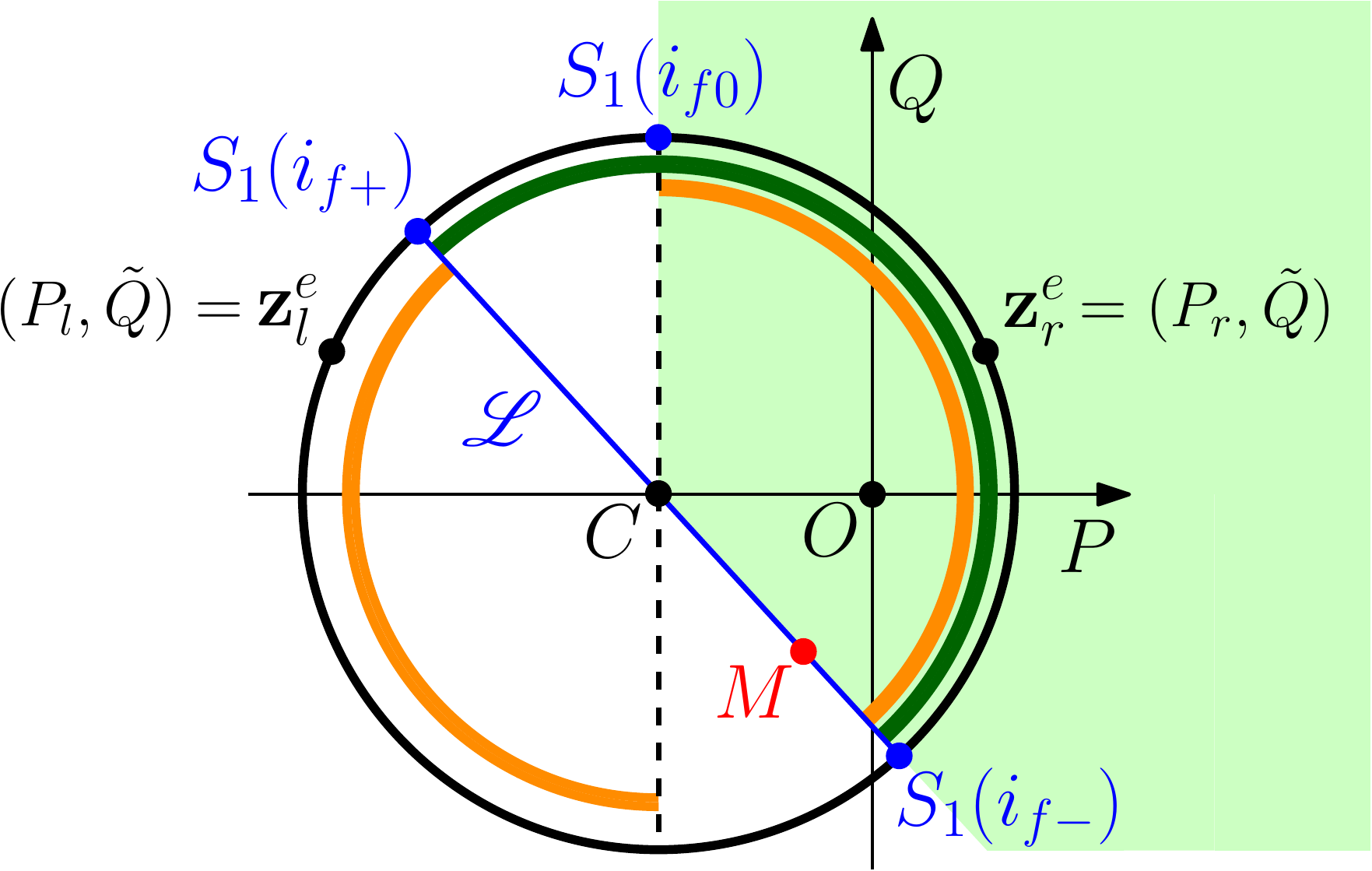} \caption{The circle
on which the equilibrium points $\zzz_r^e$ and $\zzz_l^e$ of \rfb{eq:ss5_sat}
are located for a fixed $\tilde{T}_m$. (Here we assume that $\o_g L>R$, so
that $\Lscr$ has a negative slope. The case $\o_g L\leq R$ is similar, but
derived according to Fig.~\ref{alternative}.) The green semicircle is the
stability region of \rfb{eq:ss4} on the circle, the orange arcs indicate the
region on the circle where $G'(i_f)>0$, while the light green area denotes
the stability region of \rfb{eq:ss5_sat}, for varying $P_{\rm set},\m Q_{\rm
set}$.} \label{fig:circl_expl} \end{figure}

\vspace{1mm} We now illustrate how to derive the region of stability of the
fifth order model \rfb{eq:ss5_sat} in the power plane. We assume that the
inverter parameters, as well as $V$ and $\o_g$, are known and fixed, but $P_{\rm
set}$ and $Q_{\rm set}$ can vary. Recall the notation of Theorems
\ref{circle_more}, \ref{Fakhrizadeh}. Then the coordinates
of $\zzz_r^e=S_1(i_{fr}^e)$ can be obtained from $(P_{\rm set},Q_{\rm set})$ as
follows: \begin{itemize} \item If $\o_n=\o_g$ and $v_{\rm
set}=\sqrt{\frac{2}{3}}V$, i.e., the grid is in nominal conditions, then
$S_1(i_{fr}^e)=(P_{\rm set}, Q_{\rm set})$. \item If the grid is not in nominal
conditions, $T_m$ is computed from \rfb{eq:T_m}, $\tilde{T}_m$ is given by
\rfb{eq:Tm_tilde}, and $\tilde{Q}$ is computed according to \rfb{eq:Q_tilde}.
Finally, $P$ is the larger of the two solutions of \rfb{eq:Tm_wg} and
$S_1(i_{fr}^e)=(P,\tilde{Q})$. \end{itemize} According to Proposition
\ref{prop:3.1}, $Q=\tilde{Q}$ at both the equilibrium points $\zzz_r^e$ and
$\zzz_l^e$ of \rfb{eq:ss5_sat}, and they both satisfy \rfb{eq:Tm_wg}. Hence,
$\zzz_r^e$ and $\zzz_l^e$ are located on the circle with radius $r$ and centre
$C$ given by \rfb{eq:r_C}, as show in Fig.~\ref{fig:circl_expl}.

\begin{figure*}[h!] 
\vspace{-4mm} \centering {\subfigure[Example \ref{ex:a}.]
{\includegraphics[height=0.35\textwidth]{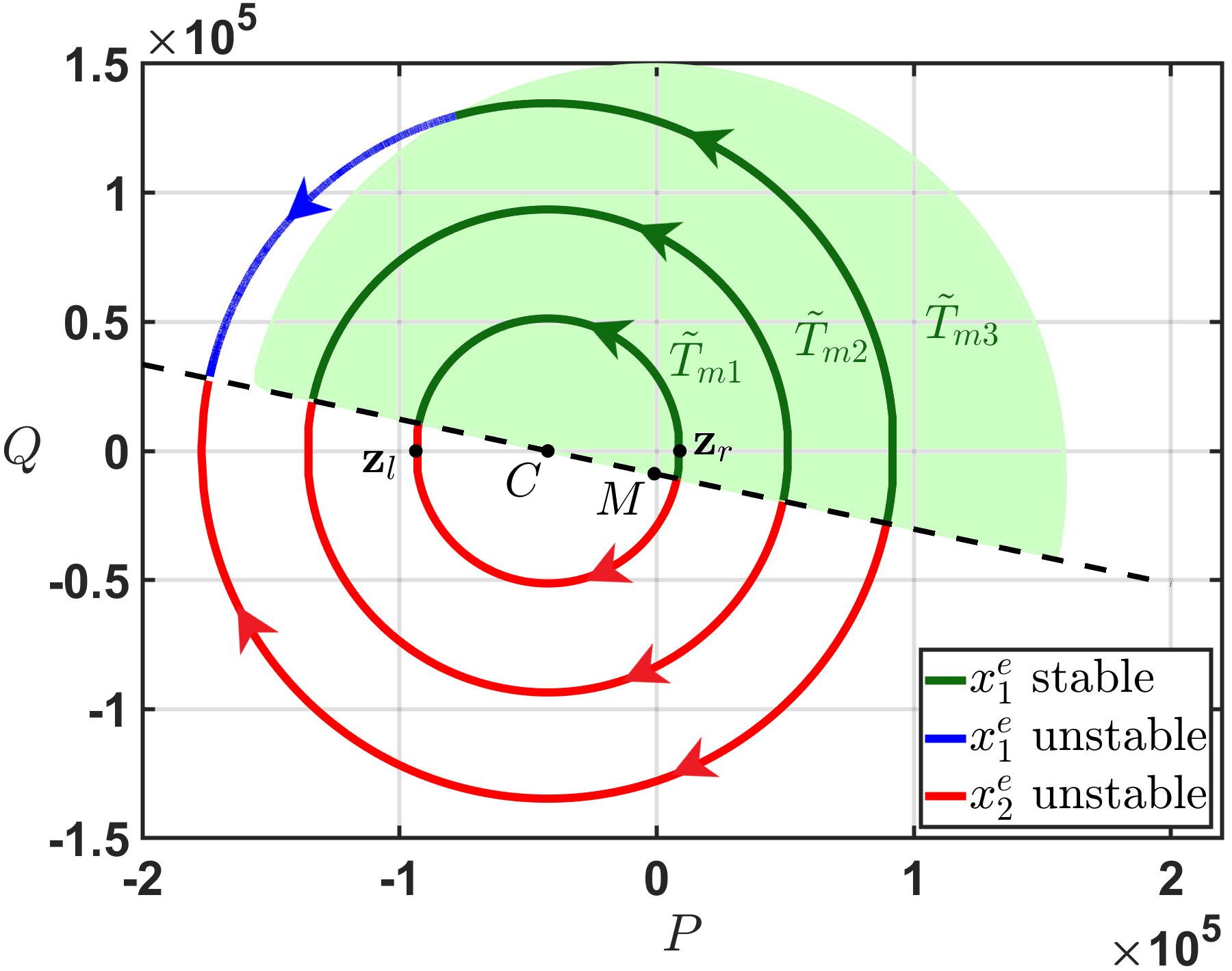}
\label{fig:ex1_circle}} \hfill \subfigure[Example
\ref{ex:b}.] {\includegraphics[height=0.35\textwidth]{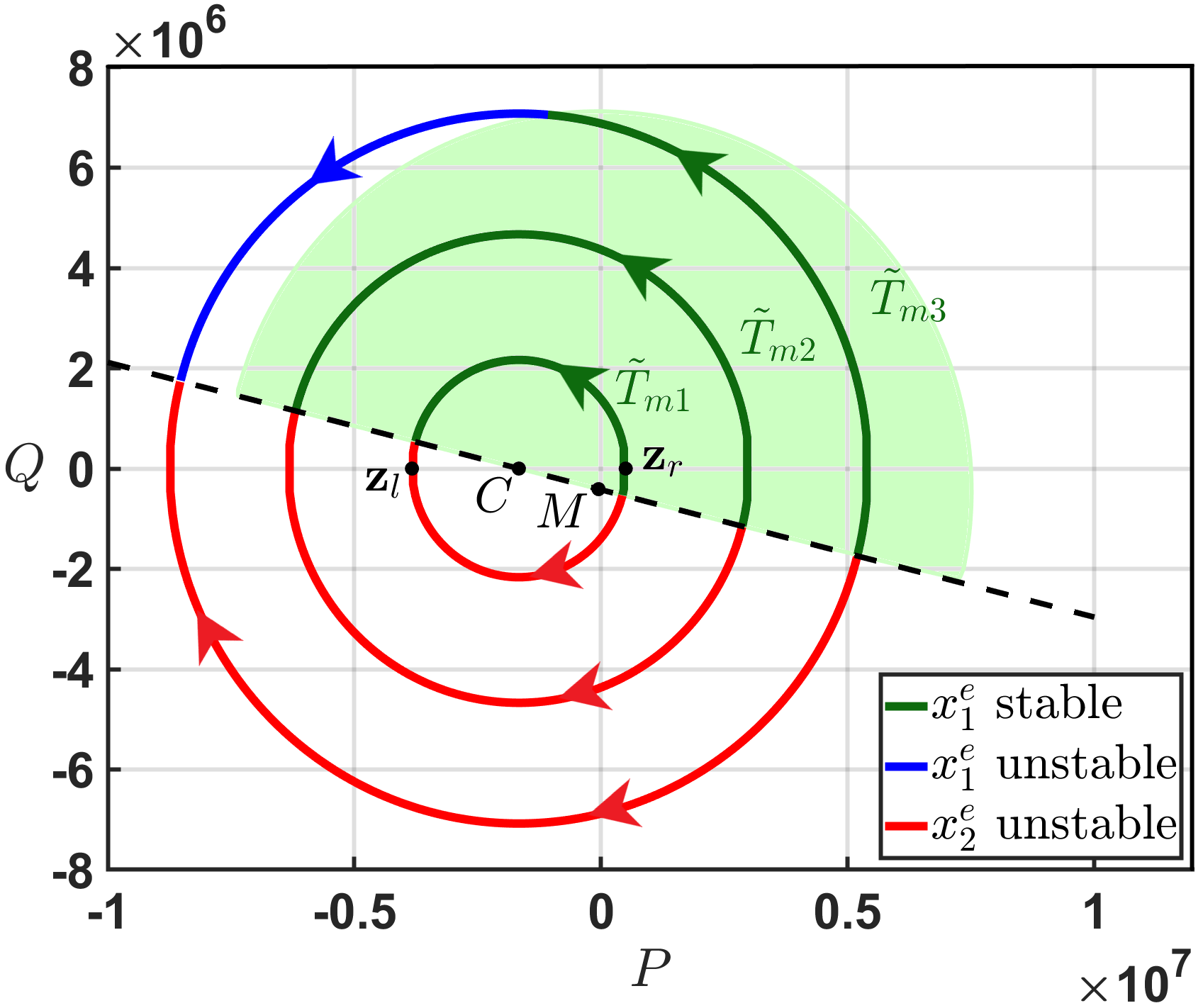}
\label{fig:ex2_circle}}} \caption{The circles are representing the
equilibrium points of the fourth order model \rfb{eq:ss4} in the power plane,
with the parameters from Example \ref{ex:a} (Subfig. a)
and from Example \ref{ex:b} (Subfig. b), for different values of
$\tilde{T}_m$ from \rfb{eq:Tm_tilde}. The legend refers to the coloured arcs,
and $\xxx_1^e$, $\xxx_2^e$ are from Proposition \ref{prop:Xi}. The light
green area denotes the maximal stability region of \rfb{eq:ss4}, with the
parameters from Example \ref{ex:a} (Subfig. a) and from
Example \ref{ex:b} (Subfig. b). The points $C,M$ are from Proposition
\ref{circle}, while $\zzz_l,\zzz_r$ are from Theorem \ref{Fakhrizadeh}.}
\end{figure*}

According to our experience (see Examples \ref{ex:a}, \ref{ex:b}), for usual
synchronverter parameters and normal operating conditions, the equilibrium
points $\xxx_1^e(i_f)$ of the fourth order model \rfb{eq:ss4} are stable for all
$i_f\in I_f$. This semicircle is indicated in dark green in
Fig.~\ref{fig:circl_expl}. On the other hand, Theorem \ref{thm:5.1} tells us
that the equilibrium points $\zzz_r^e$ of the fifth order model \rfb{eq:ss5_sat}
given by $\zzz_r^e=(\Xi(i_{fr}^e),i_{fr}^e)$, where $i_{fr}^e\in I_f$ is such
that $\xxx_1^e(i_{fr}^e)$ is stable and $G'(i_{fr}^e)>0$, are stable. Thus, if
we indicate in orange the region of the circle where $G'(i_{fr}^e)>0$ (see
Fig.~\ref{fig:circl_expl}), its intersection with the green semicircle gives the
region where the assumptions of Theorem \ref{thm:5.1} hold. From here, it
follows that, for different values of $P_{\rm set}$ and $Q_{\rm set}$ (i.e.
different values of $\tilde{T}_m$), the stability region of the resulting fifth
order model \rfb{eq:ss5_sat} is contained in the green conic sector in
Fig.~\ref{fig:circl_expl}. As will be illustrated in
Examples \ref{ex:a}, \ref{ex:b}, the stability region of \rfb{eq:ss5_sat} in the
power plane depends on the value of $\tilde{K}$. Indeed, from our computations
we see that (for fixed synchronverter parameters) the stability region is
changing for different values of $\tilde{K}$. Surprisingly, it seems that, even
though the overall stability region area $A(\tilde{K})$ is increasing for
increasing values of $\tilde{K}$, it is not true that if
$\tilde{K}_1>\tilde{K}_2$ then $A(\tilde{K}_2) \subset A(\tilde{K}_1)$.
Moreover, Theorem \ref{thm:5.1} states that if \rfb{eq:ss4} is stable and
$i_f\in I_f^+$, then also \rfb{eq:ss5_sat} must be stable for sufficiently large
values of $\tilde{K}$. However, the converse is not true. Indeed, it can happen
that \rfb{eq:ss5_sat} is stable for some values of $\tilde{K}$ in regions of the
power plane where \rfb{eq:ss4} is not, as discussed in the numerical examples of
Sect. \ref{sec6}. \vspace{10mm}


\section{Numerical Examples} \label{sec6} 

In this section, we use two examples from the
synchronverter literature to illustrate our theoretical derivations: Example
\ref{ex:a} is taken from \cite{Kust2018}, and Example \ref{ex:b} from
\cite{Natarajan2018}. The focus is the stability analysis of the fourth order
model \rfb{eq:ss4}, and of the fifth order model \rfb{eq:ss5_sat}, for varying
values of $P_{\rm set}$ and $Q_{\rm set}$. We will show how the novel
geometrical representation from Fig.~\ref{fig:circle}, \ref{alternative} is
indeed appearing naturally when studying the stability of the equilibrium points
of \rfb{eq:ss4} for $i_f\in I_f$, and we will show how the green conic sector
from Fig.~\ref{fig:circl_expl}, corresponding to the stability region of
\rfb{eq:ss5_sat}, depends on the value of $\tilde{K}$.

\vspace{-3mm}
\subsection{Low-voltage synchronverter}\label{ex:a}

We use the parameters of a synchronverter designed to supply a nominal active
power of $9$\m kW to a grid with frequency $\o_g=100\pi$\m rad/sec (50\m Hz) and
line voltage $V=230\sqrt{3}$\m Volts. This is based on a real inverter that we
have built, see \cite{Kust2018}. The parameters are: $J=0.2$\m
Kg$\cdot$m$^2$/rad, $D_p=3$\m N$\cdot$m/ (rad/sec), $L_s=2.27$\m mH,
$R_s=0.075$\m $\Om$, $K=5$\m kA, $n=25$, $D_q=0$\m VAr/Volt, $m=3.5$\m H. For
simplicity we let $v_{\rm set}= \sqrt{\frac {2}{3}}V=325.26$\m Volt, $Q_{\rm
set}=0$\m VAr, so that $\tilde Q=0$. We take $T_m=31.69$\m Nm (according to
\rfb{eq:T_m}, this mechanical torque corresponds to $P_{\rm set}=9$\m kW and
$Q_{\rm set}=0$\m VAr). We have $R=nR_s=1.875$\m $\Om$, $L=nL_s= 56.75$\m mH,
and $\phi=83.99^{\circ}$.

From Theorem \ref{Fakhrizadeh} we know that there are four equilibrium points.
We are interested in $\zzz_r^e$, $\zzz_l^e$, i.e., those corresponding to
positive $i_f$ values at the equilibrium. These can be computed as explained in
Sect. \ref{sec4}, yielding: $$ \zzz_r^e\m=\m\bbm{i_{dr}^e\\ i_{qr}^e\\ \o_g\\
\delta_r\\ i_{fr}^e} \m=\m \bbm{-15.24\\ -16.68\\ 314.16\\ 42.42^{\circ}\\
0.54}, \quad \zzz_l^e\m=\m\bbm{i_{dl}^e\\ i_{ql}^e\\ \o_g\\ \delta_l^e\\
i_{fl}^e} \m=\m \bbm{-235.04\\ -2.38\\ 314.16\\ -90.58^{\circ}\\ 3.81}.$$

We mention that if we compute the active power $P$ at the above two equilibrium
points according to \rfb{eq:P}, we get that $P_r=9$\m kW at the stable
equilibrium point (which is exactly $P_{\rm set}$), and $P_l=-93.64$\m kW at the
unstable equilibrium point. This corresponds to what we expected, based on
Theorem \ref{Fakhrizadeh}.

The equilibrium points $\zzz_r^e$ corresponding to $(P_r,\tilde{Q})$ and
$\zzz_l^e$ corresponding to $(P_l,\tilde{Q})$ are depicted in Fig.
\ref{fig:ex1_circle}, on the smallest circle, which corresponds to
$\tilde{T}_{m1}=31.69$\m N$\cdot$m, i.e., to $P_{\rm set}=9$\m kW and $Q_{\rm
set}=0$\m VAr. For this circle, we get $I_{f1}=[0.37,3.83]$\m A. In the same
figure, we also show two other circles, corresponding to the equilibrium points
of \rfb{eq:ss4} for $\tilde{T}_{m2}=261.64 $\m N$\cdot$m (i.e., $P_{\rm set\m
2}=50$\m kW and $Q_{\rm set\m 2}= 15$\m kVAr) and $\tilde {T}_{m3}=614.60$\m
N$\cdot$m (i.e., $P_{\rm set\m 3}=90$\m kW and $Q_{\rm set\m 3}=25$\m kVAr), for
which, respectively, we get $I_{f2}= [2.10, 5.56]$\m A and $I_{f3}=[3.78,
7.24]$\m A. Note that $\tilde{T}_m=T_m$, since $\o_g =\o_n$. As we can see from
Fig.~\ref{fig:ex1_circle}, while the equilibrium points $\xxx_2^e$ are always
unstable, which is a known fact according to Proposition \ref{Rouhani}, the
equilibrium points $\xxx_1^e$ in this example are always stable for reasonable
(i.e., not too large) $P_{\rm set}$ and $Q_{\rm set}$ values. This can be
checked by computing the eigenvalues of the linearizations.

The light green area in Fig.~\ref{fig:ex1_circle} indicates the stability region
of the fourth order model \rfb{eq:ss4}, which indeed covers all the relevant
$(P,Q)$ values. We mention an interesting observation: it seems from our
numerical results that the point $M$ coincides with the centre of the green
semidisk in Fig.~\ref{fig:ex1_circle}, indicating the stability region of
\rfb{eq:ss4}.

In Fig.~\ref{fig:ex1_5th} we show how the contour of the fifth order model
\rfb{eq:ss5_sat} stability region varies for different values of $\tilde{K}$. 
We use the following values: $\tilde{K}_1=2.5$\m kA$\cdot$H,
$\tilde{K}_2=14.3$\m kA$\cdot$H, $\tilde{K}_3=40$\m kA$\cdot$H, and
$\tilde{K}_4=1000$\m kA$\cdot$H. Note that $\tilde{K} _2$ is the value
corresponding to $K=5$\m kA, i.e., the one used above for the computation of
$\zzz_l^e$ and $\zzz_r^e$. Even though the overall
stability region area $A(\tilde{K})$ is increasing for increasing values of
$\tilde{K}$, it is not true that if $\tilde{K}_1>\tilde{K}_2$, then
$A(\tilde{K}_2)\subset A(\tilde{K}_1)$, as is clear from
Fig.~\ref{fig:ex1_5th}. We mention that, for $\tilde{K}\to\infty$, it seems
from our numerical results that the region of stability of \rfb{eq:ss5_sat}
coincides with the intersection of the green sector from
Fig.~\ref{fig:circl_expl} and of the stability region of \rfb{eq:ss4}. This can
be observed in Fig.~\ref{fig:ex1_5th}, where, for increasing values of
$\tilde{K}$, the stability region contours approach the boundary of the light
green area.

\vspace{-2mm}
\subsection{High-voltage synchronverter}\label{ex:b}

We consider a synchronverter from \cite{Natarajan2018} that supplies a nominal
active power of $500$\m kW to a grid with frequency $\o_g=100\pi$\m rad/sec
(50\m Hz) and line voltage $V=6000\sqrt{3}$\m Volts. The parameters are:
$J=20.26\,$\m Kg$\cdot $m$^2$/rad, $D_p=168.87$\m N$\cdot$m/(rad/sec),
$L_s=27.5$\m mH, $R_s= 1.08$\m $\Om$, $K=5000$\m A, $n=30$, $D_q=0$\m VAr/Volt,
$m=33$\m H. As previously, we let $v_{\rm set}=\sqrt{\frac{2}{3}}V=8485.3$\m
Volt, $Q_{\rm set}=0$\m VAr, so that $\tilde Q=0$. The mechanical torque
$T_m=1.83$\m kN$\cdot$m (according to \rfb{eq:T_m}) corresponds to $P_{\rm
set}=500$\m kW and $Q_{\rm set}=0$\m VAr. We have $R=nR_s= 32.4$\m $\Om$,
$L=nL_s=825$\m mH, and $\phi= 82.87^{\circ}$.

The two equilibrium points with positive $i_f$ values are: $$
\zzz_r^e\m=\m\bbm{i_{dr}^e\\ i_{qr}^e\\ \o_g\\ \delta_r\\ i_{fr}^e} \m=\m
\bbm{-34.73\\ -33.29\\ 314.16\\ 46.21^{\circ}\\ 1.67}, \quad
\zzz_l^e\m=\m\bbm{i_{dl}^e\\ i_{ql}^e\\ \o_g\\ \delta_l^e\\ i_{fl}^e} \m=\m
\bbm{-368.81\\ -6.01\\ 314.16\\ -90.93^{\circ}\\ 9.22}.$$

Again, if we compute the active power $P$ at the above two equilibrium points
according to \rfb{eq:P}, we get that $P_r=500$\m kW at the stable equilibrium
point (which is exactly $P_{\rm set}$), and $P_l=-3.83$\m MW at the unstable
equilibrium point. In the following, we perform the same stability analysis of
Example \ref{ex:a}.

The equilibrium points $\zzz_r^e$ corresponding to $(P_r,\tilde{Q})$ and
$\zzz_l^e$ corresponding to $(P_l,\tilde{Q})$ are shown in
Fig.~\ref{fig:ex2_circle}, on the smallest circle, which corresponds to
$\tilde{T}_{m1}=1.83$\m kN$\cdot$m, i.e., to $P_{\rm set}=500$\m kW and $Q_{\rm
set}=0$\m VAr. For this circle, we get $I_{f1}=[1.21, 9.29]$\m A. In the same
figure, we also represent two other circles, corresponding to the equilibrium
points of \rfb{eq:ss4} for $\tilde{T}_{m2}=18.18$\m kN$\cdot$m (i.e., $P_{\rm
set\m 2}=3000$\m kW and $Q_{\rm set\m 2}=200$\m kVAr and
$\tilde{T}_{m3}=45.19$\m kN$\cdot$m (i.e., $P_{\rm set\m 3}=5400$\m kW and
$Q_{\rm set\m 3}=400$\m kVAr), for which, respectively, we get $I_{f2}=[7.28,
15.36]$\m A and $I_{f3}=[13.12, 21.20]$\m A. Also with these synchronverter
values, it is clear from Fig.~\ref{fig:ex2_circle} that the equilibrium points
$\xxx_1^e$ are always stable for reasonable $P_{\rm set}$ and $Q_{\rm set}$
values. (Only for $\tilde{T}_{m3}$ can we see a blue arc appearing.) This is
confirmed by the light green area in Fig.~\ref{fig:ex2_circle}, indicating the
stability region of the fourth order model \rfb{eq:ss4}. Also in this case, the
point $M$ coincides with the centre of the green semidisk in Fig.
\ref{fig:ex2_circle}, indicating the stability region of \rfb{eq:ss4}.

\begin{figure}[h] 
\centering {\subfigure[Example \ref{ex:a}]
{\includegraphics[height=0.30\textwidth]{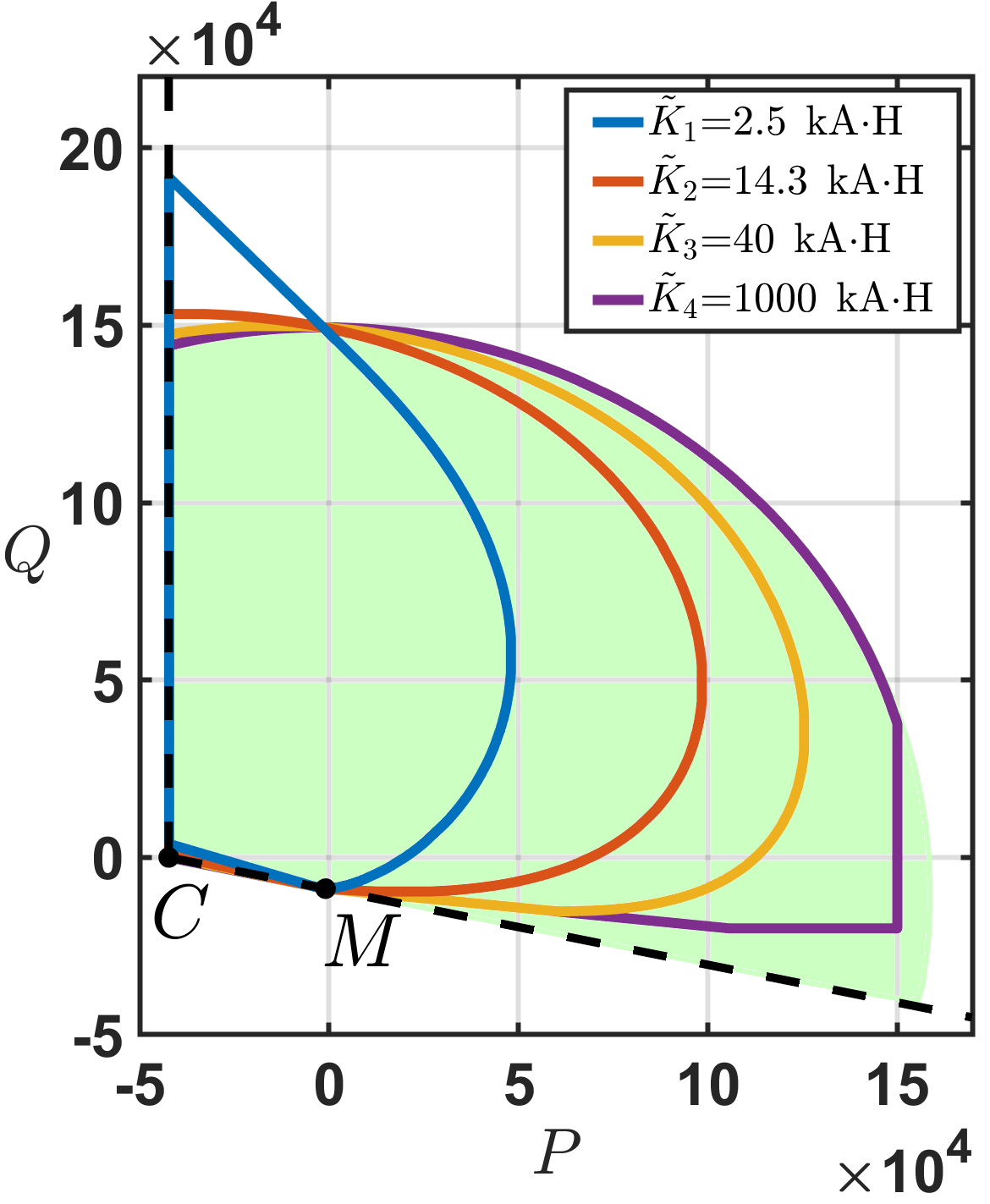}
\label{fig:ex1_5th}} \hfil \subfigure[Example
\ref{ex:b}] {\includegraphics[height=0.30\textwidth]{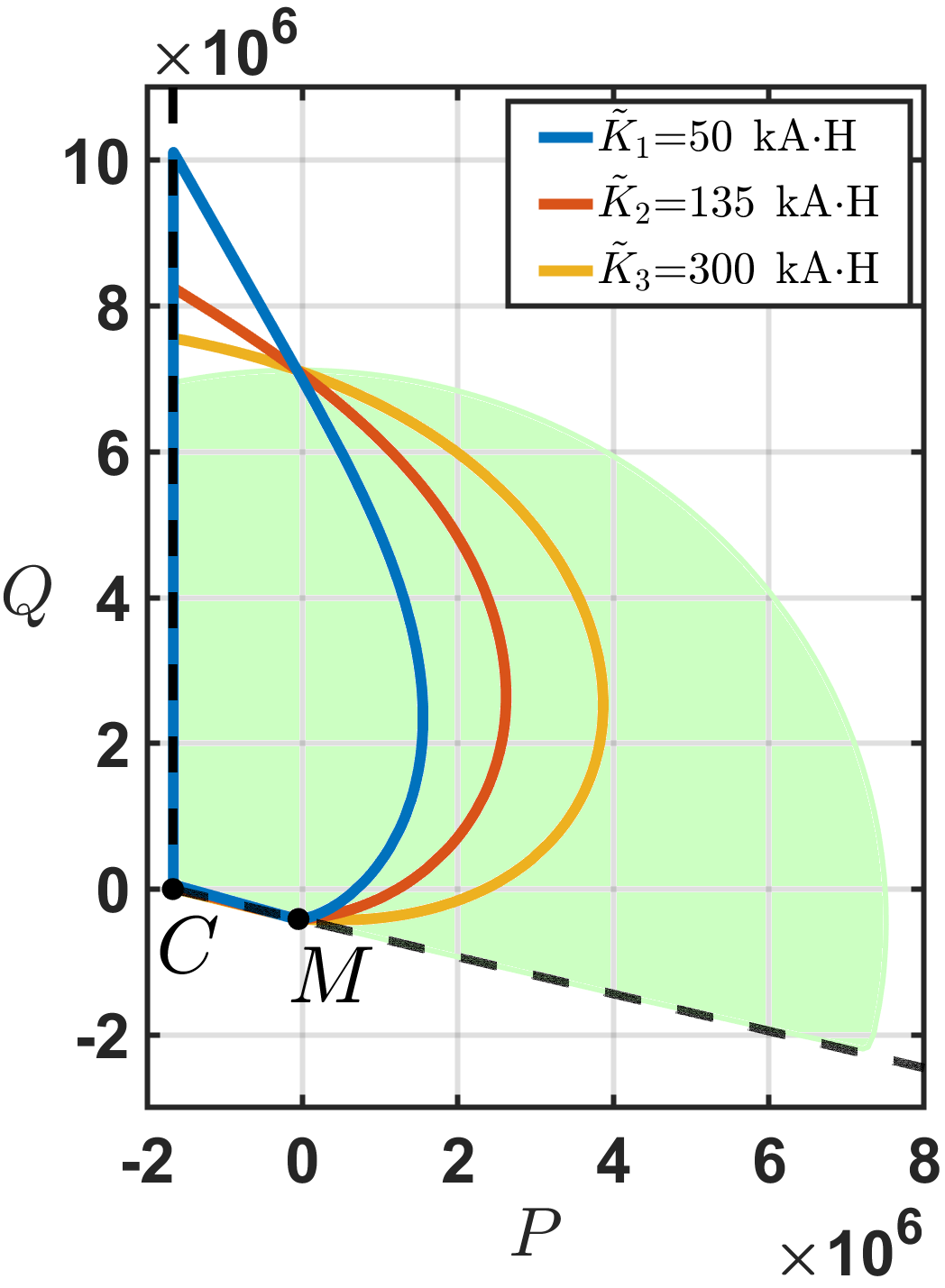}
\label{fig:ex2_5th}}} \caption{The coloured shapes correspond to the contours
of the stability region of the fifth order model \rfb{eq:ss5_sat}, with the
parameters from Example \ref{ex:a} (Subfig.~a) and from
Example \ref{ex:b} (Subfig.~b), for different values of $\tilde{K}$. The
points $C,M$ are from Proposition \ref{circle}, while the light green area
indicates the intersection between the green sector from
Fig.~\ref{fig:circl_expl} and the maximal stability region of \rfb{eq:ss4},
with the parameters from Example \ref{ex:a} (Subfig.~a) 
and from Example \ref{ex:b} (Subfig.~b).} \vspace{-3mm} \end{figure}

In Fig.~\ref{fig:ex2_5th} we show how the contours of the fifth order model
\rfb{eq:ss5_sat} stability region vary for different values of $\tilde{K}$. We
use the following values: $\tilde{K}_1=50$\m kA$\cdot$H, $\tilde{K}_2=135$\m
kA$\cdot$H, and $\tilde{K}_3=300$\m kA$\cdot$H. Note that
$\tilde{K}=\tilde{K}_2$ is the value corresponding to $K=5$\m kA, i.e., the one
used above for the computation of $\zzz_l^e$ and $\zzz_r^e$. Also in this case,
it is not true that if $\tilde{K}_1>\tilde{K}_2$, then $A(\tilde{K}_2)\subset
A(\tilde{K}_1)$, as is clear from Fig.~\ref{fig:ex2_5th}. Moreover, we observe,
again, that for $\tilde{K}\to\infty$ the contours are approaching the green
light area, indicating the intersection between the green sector from
Fig.~\ref{fig:circl_expl} and the stability region of \rfb{eq:ss4}.

\vspace{-3mm}
\section{Conclusions} \label{sec7} 

We have formulated a fifth order model for a grid-connected synchronverter, when
the grid is considered to be an infinite bus. Conditions ensuring the existence
of its equilibrium points have been derived, and a novel geometrical
representation has been introduced. This representations links the region of
stability of the fourth order model from \cite{Natarajan2018,Natarajan2017},
with the region of stability of our fifth order model. Moreover, using singular
perturbation methods, we have derived sufficient conditions guaranteeing the
existence of (local) exponentially stable equilibrium points for the fifth order
model. Finally, the validity of our theoretical results has been proved using
two numerical examples coming from the synchronverter literature.

 \vspace{5mm}

\small{ {\bf Pietro Lorenzetti} is an Early Stage Researcher within the Marie
Curie ITN project ``ConFlex'', who focuses his research on nonlinear control. 
Pietro has completed the bachelor degree in ``Computer engineering and
automation'' at Universita Politecnica delle Marche, in Ancona. In 2015 he
graduated with honors and he moved to Torino, where he enrolled the master
degree in ``Mechatronic Engineering'' at Politecnico di Torino. In the same
year, he also joined the double-degree program ``Alta Scuola Politecnica'', a
highly selective joined program between Politecnico di Torino and Politecnico di
Milano. In 2017 he graduated in both Politecnico di Milano and Politecnico di
Torino, with honors. His research interests include nonlinear systems, nonlinear
control, and power system stability.}

\small{ {\bf Zeev Kustanovich} received the B.Sc. degree from Ben Gurion
University of the Negev, Beer Sheva, Israel, in 1997, and the M.Sc. degree from
Technion, Haifa, Israel, in 2003 both in electrical engineering. Since 2003 he
is Senior Electrical Engineer at the Israel Electricity Company. In 2018, he
started his Ph.D. with the Power Electronics for Renewable Energy group in Tel
Aviv University, Israel. His main research interests include power systems,
renewable energy, control theory and applicatons to power system stability.}

\small{ {\bf Shivprasad Shivratri} received the B.Sc. degree in electrical and
electronics engineering from Tel Aviv University, Israel, in 2018. In 2018, he
started his M.Sc. with the Power Electronics for Renewable Energy group in Tel
Aviv University and he graduated in May 2021. His research interests include
control techniques in power systems and control theory.}

\small{ {\bf George Weiss} received the MEng degree in control engineering from
the Polytechnic Institute of Bucharest, Romania, in 1981, and the Ph.D. degree
in applied mathematics from the Weizmann Institute, Rehovot, Israel, in 1989. He
was with Brown University, Providence, RI, Virginia Tech, Blacksburg, VA,
Ben-Gurion University, Beer Sheva, Israel, the University of Exeter, U.K., and
Imperial College London, U.K. His current research interests include distributed
parameter systems, operator semigroups, passive and conservative systems (linear
and nonlinear), power electronics, microgrids, repetitive control, sampled data
systems, and wind-driven power generators. He is leading research projects for
the European Commission and for the Israeli Ministry of Infrastructure, Energy
and Water.} \end{document}